\documentstyle[epsfig,amsfonts,amscd]{elsart}

\input{xypic}

\DeclareSymbolFont{EulerScript}{U}{eus}{m}{n}
\DeclareSymbolFontAlphabet\mathscr{EulerScript}

\newcommand{\rk}{{\hbox{\rm rk}}}
\newcommand{\cok}{{\hbox{\rm coker}}}
\newcommand{\Ann}{{\hbox{\rm Ann}}}
\newcommand{\diag}{{\hbox{\rm diag}}}
\newcommand{\mod}{{\,\hbox{\rm mod}\,}}
\newcommand{\Z}{{\Zset}}
\newcommand{\N}{{\mathscr N}}
\newcommand{\R}{{\Rset}}
\newcommand{\im}{{\hbox{Im}}}
\newcommand{\Hom}{{\hbox{Hom}}}
\newcommand{\eps}{\varepsilon}

\begin{document}
	
\begin{frontmatter}

\title{A Lagrangian representation of tangles}

\author{David Cimasoni\thanksref{SNSF}}
\address{Section de Math\'ematiques, Universit\'e de Gen\`eve, 2--4 
rue du Li\`evre, 1211 Gen\`eve 24, Switzerland}
\author{Vladimir Turaev\thanksref{CA}}
\address{IRMA, UMR 7501 CNRS/ULP, 7 rue Ren\'e Descartes, 67084 
Strasbourg Cedex, France}
\thanks[SNSF]{Supported by the Swiss National Science Foundation.}
\thanks[CA]{Corresponding author. Tel.:+33 (0)3 90 240 149;
fax: +33 (0)3 90 240 328.\\
{\em E-mail address:\/} turaev@math.u-strasbg.fr}

\begin{abstract}
We construct a functor from the category of oriented tangles in $\R^3$ 
to the category of Hermitian modules and Lagrangian
relations over $\Z[t,t^{-1}]$. This functor extends the Burau 
representations of the braid groups and its generalization
to string links due to Le Dimet. 
\end {abstract}


\end{frontmatter}

\section{Introduction}

The aim of this paper is to generalize the classical Burau 
representation of braid groups to tangles.
The Burau representation is a homomorphism from the group of braids on 
$n$ strands to the group of $(n\times n)$-matrices
over the ring $\Lambda=\Z[t,t^{-1}]$, where $n$ is a positive integer. 
This representation has been extensively
studied by various authors since the foundational work of Burau 
\cite{Bur}. In the last 15 years, new important
representations of braid groups came to light, specifically those 
associated with the Jones knot polynomial,
$R$-matrices, and ribbon categories. These latter representations do 
extend to tangles, so it is natural to ask whether
the Burau representation has a similar property.

An extension of the Burau representation to a certain class of tangles 
was first pointed out by Le Dimet \cite{LeD}.
He considered so-called `string links', which are tangles whose all 
components are intervals going from the bottom
to the top but not necessarily monotonically. The string links on $n$ 
strands form a monoid with respect to the usual
composition of tangles. Le Dimet's work yields a homomorphism of this 
monoid into the group of $(n\times n)$-matrices
over the quotient field of $\Lambda$. For braids, this gives the Burau 
representation. The construction of Le Dimet also
applies to colored string links, giving a generalization of the 
Gassner representation of the pure braid group.
These representations of Le Dimet were studied by Kirk, Livingston and 
Wang \cite{KLW} (see also \cite{LTW,SW}).

To extend the Burau representation to arbitrary oriented tangles, we 
first observe that oriented tangles do not form a group
or a monoid but rather a category ${\mathbf Tangles}$ whose objects are 
finite sequences of $\pm 1$. An extension of the Burau
representation to ${\mathbf Tangles}$ should be a functor from 
${\mathbf Tangles}$ to some algebraically defined category.
We show that the relevant algebraic category is the one of Hermitian 
$\Lambda$-modules and Lagrangian relations.
Our principal result is a construction of a functor from 
${\mathbf Tangles}$ to this category. For braids and string links,
our constructions are equivalent to those of Burau and Le Dimet.

The appearance of Lagrangian relations rather than homomorphisms is 
parallel to the following 
well-known observations concerning cobordisms. Generally speaking, a 
cobordism $(W,M_-,M_+)$ does not induce a homomorphism
from the homology (with any coefficients) of the  bottom base $M_-$ to 
the homology of  the top base $M_+$.
However, the kernel of the inclusion homomorphism $H_\ast(M_-)\oplus 
H_\ast(M_+) \to H_\ast(W)$  can be viewed as a morphism
from $H_\ast(M_-)$ to $ H_\ast(M_+)$ determined by $W$. This  kernel 
is Lagrangian with respect to the usual 
intersection form in homology. These observations suggest a definition 
of a Lagrangian category over any integral domain
with involution. Applying these ideas to the infinite cyclic covering 
of the tangle exterior, we obtain our functor from the
category of tangles to the category of Lagrangian relations over 
$\Lambda$.

Note that recently, a most interesting representation of braid groups 
due to
R. Lawrence was shown to be faithful by S. Bigelow and D. Krammer.
We do not know whether this representation extends to tangles.

The organization of the paper is as follows. In Section 2, we 
introduce the category ${\mathbf Lagr_\Lambda}$ of
Lagrangian relations over the ring $\Lambda$. In Section 3, we define 
our functor ${\mathbf Tangles}\to {\mathbf Lagr_\Lambda}$.
Section 4 deals with the proof of three technical lemmas stated in the 
previous section. In Section 5, we discuss the case of
braids and string links. Section 6 concerns technical questions about 
the Lagrangian relations associated with tangles.
In Section 7, we discuss connexions between these Lagrangian relations 
and the Alexander polynomial of the link obtained
as the closure of the tangle. (These connexions are traditionally 
studied in this context.) Finally, Section 8 outlines a
multivariable generalization of the theory as well as a 
high-dimensional version.

\section{Category of Lagrangian relations}

Fix throughout this section an integral domain $\Lambda$ (i.e., a 
commutative ring with unit and
without zero-divisors) with ring involution 
$\Lambda\to\Lambda,\lambda\mapsto\tilde{\lambda}$.

\subsection{Hermitian modules}

A {\em skew-hermitian form} on a $\Lambda$-module $H$ is a form 
$\omega\colon H\times H\to\Lambda$ such that for all
$x,x',y\in H$ and all $\lambda,\lambda'\in\Lambda$,
\begin{enumerate}
\item{$\omega(\lambda 
x+\lambda'x',y)=\lambda\omega(x,y)+\lambda'\omega(x',y),$}
\item{$\omega(x,y)=-\widetilde{\omega(y,x)}$.}
\end{enumerate}
Such a form is called {\em non-degenerate} when it satisfies:
\begin{enumerate}\setcounter{enumi}{2}
\item{If $\omega(x,y)=0$ for all $y\in H$, then $x=0$.}
\end{enumerate}
A {\em Hermitian $\Lambda$-module} is a finitely generated 
$\Lambda$-module
$H$ endowed with a non-degenerate skew-hermitian form $\omega$.
The same module $H$ with the opposite form $-\omega$ will be denoted 
by $-H$. 
Note that a Hermitian $\Lambda$-module is always torsion-free.

For a submodule $A\subset H$, denote by $\Ann(A)$ the annihilator of 
$A$ with respect to $\omega$, that is, the
module $\{x\in H \mid \omega(x,a)=0 \hbox{ for all $a\in A$}\}$. Set 
\[
\overline A=\{x\in H \mid \lambda x\in A \hbox{ for a non-zero 
$\lambda\in\Lambda$}\}.
\]
Clearly $A\subset \overline A$ and 
$\overline{\Ann(A)}=\Ann(A)=\Ann(\overline{A})$. 

We say that a submodule $A$ of $H$ is {\em isotropic} if $A\subset\Ann 
(A)$.
Observe that then $A\subset \overline A\subset\Ann (A)$. 
We say that a submodule $A$ of $H$ 
is {\em Lagrangian} if $\overline{A}=\Ann(A)$. Note that $A$ is 
Lagrangian if and
only if $\overline A$ is Lagrangian.

\begin{lem}\label{lemma2}
For any submodule $A$ of a Hermitian $\Lambda$-module $H$,
\[
\Ann(\Ann(A))=\overline A.
\]
\end{lem}
\begin{pf}
Let $Q=Q(\Lambda)$ denote  the field of fractions of $\Lambda$. Given 
a
$\Lambda$-module $F$, denote by $F_Q$ the vector space 
$F\otimes_\Lambda Q$.
Note that the kernel of the natural homomorphism $F\to F_Q$ is the
$\Lambda$-torsion $Tors_\Lambda F\subset F$. 

The form  $\omega$  uniquely  extends
to a skew-hermitian form $H_Q\times H_Q\to Q$. 
Given a linear subspace $V$ of $H_Q$, let $\Ann_Q(V)$ be the 
annihilator of $V$
with respect to the latter form. Observe that
$\Ann_Q(\Ann_Q(V))=V$.
 Indeed, one inclusion is trivial and the other one 
 follows from dimension count, since 
$\dim(\Ann_Q(V))=\dim(H_Q)-\dim(V)$.
 
The inclusion   $A\hookrightarrow H$ induces an inclusion  $A_Q 
\hookrightarrow
H_Q$.
Since $H$ is torsion-free,   
 $H\subset H_Q$ (and $A\subset A_Q$). Clearly,  $\overline A=A_Q\cap 
H$ and 
$\Ann(A)_Q=\Ann_Q(A_Q)$. Replacing in the latter formula $A$ with 
$\Ann (A)$,
we obtain
\[
\Ann(\Ann(A))_Q=\Ann_Q(\Ann(A)_Q)=\Ann_Q(\Ann_Q(A_Q))=A_Q.
\]
Therefore  
\[
\overline A=A_Q\cap H=\Ann(\Ann(A))_Q\cap
H=\overline{\Ann(\Ann(A))}=\Ann(\Ann(A)),
\]
and the lemma is proved.
\qed\end{pf}

\begin{lem}\label{lemma3}
For any submodules $A,B\subset H$, we have
\begin{eqnarray*}
\Ann(A+B)&=&\Ann(A)\cap\Ann(B),\\
\Ann(A\cap B)&=&\overline{\Ann(A)+\Ann(B)}.
\end{eqnarray*}
\end{lem}
\begin{pf}
The first equality is obvious, and implies
\begin{eqnarray*}
\Ann(\Ann(A)+\Ann(B))&=&\Ann(\Ann (A))\cap\Ann(\Ann(B))\\
	&=&\overline{A}\cap\overline{B}=\overline{A\cap B}.
\end{eqnarray*}
Therefore
\[
\Ann(A\cap B)=\Ann(\overline{A\cap B})=\Ann(\Ann(\Ann(A)+\Ann(B))),
\]
which is equal to $\overline{\Ann(A)+\Ann(B)}$ by Lemma \ref{lemma2}. 
\qed\end{pf}

\begin{lem}\label{lemma:quotient}
For any submodules $A\subset B\subset H$, we have
$\overline{B}/A=\overline{B/A}$.
\end{lem}
\begin{pf}
Consider the  canonical projection $\pi\colon H\to H/A$. 
Clearly,
\[
\pi(\overline{B})=\{\xi\in H/A \mid \lambda \xi\in B/A \hbox{ for a 
non-zero 
$\lambda\in\Lambda$}\}=\overline{B/A}.\]
Also
$\ker(\pi\vert_{\overline{B}})=\ker(\pi)\cap\overline{B}=A\cap
\overline{B}=A$. 
Hence $\overline{B}/A=\overline{B/A}$.
\qed\end{pf}

\subsection{Lagrangian contractions}

The  results above in hand, it is easy to develop the theory of 
Lagrangian contractions 
and Lagrangian relations over $\Lambda$ by mimicking the well-known 
theory over $\R$ (see, for instance, 
\cite[Section IV.3]{Tur}).

Let $(H,\omega)$ be a Hermitian $\Lambda$-module as above.
 Let $A$ be an isotropic submodule of $H$ such that $A=\overline A$.
  Denote by $H\vert A$ the quotient module $\Ann(A)/A$ with the 
skew-hermitian form
\[
(x\mod A,y\mod A)=\omega(x,y).
\]
For a submodule $L\subset H$, set
\[
L\vert A=((L+A)\cap \Ann(A))/A \,\subset H\vert A.
\] 
We say that $L\vert A$ is obtained from $L$ by {\em contraction along} 
$A$. 

\begin{lem}\label{lemma:contraction}
$H\vert A$ is a Hermitian $\Lambda$-module.
If $L$ is a Lagrangian submodule of $H$, then $L\vert A$ is a 
Lagrangian submodule of $H\vert A$.
\end{lem}
\begin{pf}
  To check that the form on $H\vert A$ is non-degenerate, pick 
$x\in\Ann(A)$ such that $ \omega(x,y)=0$
  for all $y\in\Ann(A)$. Then, $x\in \Ann(\Ann(A))=\overline A=A$ so 
that   $x\mod A=0$.

To prove the second claim of the lemma, set $B=(L+A)\cap\Ann(A)\subset 
H$. We claim that $B$ is Lagrangian.
Since both $A$ and $L$ are isotropic, it is easy to check that 
$B\subset\Ann(B)$ and therefore $\overline B
\subset\Ann(B)$. Let us verify the opposite inclusion.
  Lemmas \ref{lemma2} and \ref{lemma3} imply that 
\begin{eqnarray*}
\Ann(B)&=&\Ann((L+A)\cap\Ann(A))=\overline{\Ann(L+A)+\Ann(\Ann(A))}\\
        &\subset&\overline{\Ann(L)+\overline A}=\overline{L+A}.
\end{eqnarray*}
Since $A\subset B$, we have  $\Ann(B)\subset \Ann(A)$ and therefore
\[
\Ann(B) 
\subset\overline{L+A}\cap\Ann(A)=\overline{(L+A)\cap\Ann(A)}=\overline 
B.
\]
Thus $B$ is Lagrangian.  This implies that $\Ann(B/A)=\overline{B}/A$, 
which is equal to $\overline{B/A}$
by Lemma \ref{lemma:quotient}. So $B/A$ is Lagrangian.
\qed\end{pf}

\subsection{Categories of Lagrangian relations}

Let $H_1,H_2$ be   Hermitian $\Lambda$-modules. A {\em Lagrangian 
relation} between $H_1$ and $H_2$ is a Lagrangian
submodule of $(-H_1)\oplus H_2$ (the latter is a Hermitian 
$\Lambda$-module in the obvious way). For a Lagrangian
relation $N\subset (-H_1)\oplus H_2$, we shall use the notation 
$N\colon H_1\Rightarrow H_2$.

For a Hermitian $\Lambda$-module $H$, the submodule of $H\oplus H$
\[
\diag_H=\{ h\oplus h \in (-H)\oplus H \mid h \in H \}
\]
is clearly a Lagrangian relation $H\Rightarrow H$. It is called the 
{\em diagonal Lagrangian relation}. 
Given two Lagrangian relations $N_1\colon H_1\Rightarrow H_2$ and 
$N_2\colon H_2\Rightarrow H_3$,
their composition is the following submodule of $(-H_1)\oplus H_3$:
\[
N_2N_1=\{h_1\oplus h_3\mid \hbox{$h_1\oplus h_2 \in N_1$ and 
$h_2\oplus h_3 \in N_2$ for a certain $h_2\in H_2$}\}.
\]

\begin{lem}\label{lemma:Lag}
The composition of two Lagrangian relations is a Lagrangian relation.
\end{lem}
\begin{pf}
Given two Lagrangian relations $N_1\colon H_1\Rightarrow H_2$ and 
$N_2\colon H_2\Rightarrow H_3$, consider
the Hermitian $\Lambda$-module $H=(-H_1)\oplus H_2\oplus (-H_2)\oplus 
H_3$ and its isotropic submodule 
\[
A=0\oplus\diag_{H_2}\oplus 0=\{0\oplus h\oplus h\oplus 0\mid h\in 
H_2\}.
\]
Note that $\overline A=A$. It follows from the non-degeneracy of $H_2$ 
that $\Ann(A)=(-H_1)\oplus \diag_{H_2} \oplus H_3$.
Therefore $H\vert A=(-H_1)\oplus H_3$. Observe that $N_2N_1=(N_1\oplus 
N_2)\vert A$.
Lemma \ref{lemma:contraction} implies that ${N_2N_1}$ is a Lagrangian 
submodule of $(-H_1)\oplus H_3$.
\qed\end{pf}

\begin{thm}\label{thm:category}
Hermitian $\Lambda$-modules, as objects, and Lagrangian relations, as 
morphisms, form a category.
\end{thm} 
\begin{pf}
The composition law is well-defined by Lemma \ref{lemma:Lag} and the 
associativity follows from the definitions.
Finally, the role of the identity morphisms is played by the diagonal 
Lagrangian relations.
\qed\end{pf}

We shall call this category the {\em category of Lagrangian relations 
over $\Lambda$}. It will be denoted by
${\mathbf Lagr_\Lambda}$. Let us conclude this section with the 
definition of another Lagrangian category
${\mathbf \overline{Lagr}_\Lambda}$ closely related to the former one. 
The objects of ${\mathbf \overline{Lagr}_\Lambda}$
are Hermitian $\Lambda$-modules, and the morphisms are Lagrangian 
relations $N$ such that $\overline{N}=N$.
Finally, the composition between two morphisms $N_1\colon 
H_1\Rightarrow H_2$ and $N_2\colon H_2\Rightarrow H_3$
is defined by $N_2\circ N_1=\overline{N_2N_1}\colon H_1\Rightarrow 
H_3$.

\begin{lem}\label{lemma:comp'}
Given two submodules $N_1\subset H_1\oplus H_2$ and $N_2\subset 
H_2\oplus H_3$,
\[
\overline{N_2N_1}=\overline{\overline N_2\overline N_1}.
\]
\end{lem}
\begin{pf}
Consider an element $h_1\oplus h_3$ of $\overline N_2\overline N_1$. 
By definition, $h_1\oplus h_2\in\overline{N_1}$
and $h_2\oplus h_3\in\overline{N_2}$ for some $h_2\in H_2$, so 
$\lambda_1(h_1\oplus h_2)\in N_1$ and
$\lambda_2(h_2\oplus h_3)\in N_2$ for some $\lambda_1,\lambda_2\neq 
0$. Then  $\lambda_1\lambda_2(h_1\oplus h_3)\in N_2N_1$,
so $h_1\oplus h_3\in\overline{N_2N_1}$. Hence, $\overline N_2\overline 
N_1\subset\overline{N_2N_1}$. Taking the
closure on both sides, we get $\overline{\overline N_2\overline 
N_1}\subset\overline{N_2N_1}$. The opposite inclusion is obvious.
\qed\end{pf}

\begin{thm}
${\mathbf \overline{Lagr}_\Lambda}$ is a category, and the map 
$N\mapsto\overline{N}$ defines a functor
\[
{\mathbf Lagr_\Lambda}\stackrel{\mathscr 
J}{\longrightarrow}{\mathbf \overline{Lagr}_\Lambda}.
\]
\end{thm} 
\begin{pf}
The composition law in ${\mathbf \overline{Lagr}_\Lambda}$ is 
well-defined by Lemma \ref{lemma:Lag};
let us check that it is associative. Consider Lagrangian relations 
$N_1\colon H_1\Rightarrow H_2$, $N_2\colon H_2\Rightarrow H_3$,
and
$N_3\colon H_3\Rightarrow H_4$ such that $\overline{N_i}=N_i$ for 
$i=1,2,3$. By Lemma \ref{lemma:comp'},  
\[
N_3\circ(N_2\circ N_1) = \overline{N_3\overline{N_2N_1}} = 
\overline{\overline{N}_3\overline{N_2N_1}} = \overline{N_3(N_2N_1)}.
\]
Similarly, $(N_3\circ N_2)\circ N_1=\overline{(N_3N_2)N_1}$. The 
result now follows from the associativity of
the composition in ${\mathbf Lagr_\Lambda}$.
The role of the identity morphisms is played by the diagonal 
Lagrangian relations. Indeed, for any
Lagrangian relation $N\colon H_1\Rightarrow H_2$ such that 
$\overline{N}=N$,
\[
\diag_{H_2}\circ N=\overline{\diag_{H_2}N}=\overline N=N.
\]
Similarly, $N\circ\diag_{H_1}=N$.
Finally, let us check that the map $N\mapsto\overline{N}$ is 
functorial. Consider two Lagrangian relations
$N_1\colon H_1\Rightarrow H_2$ and $N_2\colon H_2\Rightarrow H_3$. By 
Lemma \ref{lemma:comp'},
\[
\overline{N_2}\circ\overline{N_1}=\overline{\overline N_2\overline 
N_1}=\overline{N_2N_1}.
\]
This finishes the proof.
\qed\end{pf}

\subsection{Lagrangian relations from unitary isomorphisms}

By the {\em graph} of a homomorphism $f\colon A\to B$ of abelian 
groups, we mean the set
\[
\Gamma_f=\{a\oplus f(a)\vert a\in A\}\subset A\oplus B.
\]
Let $H_1$, $H_2$ be Hermitian $\Lambda$-modules. Consider the 
Hermitian $Q$-modules $H_1\otimes Q$ and $H_2\otimes Q$, where 
$Q=Q(\Lambda)$ is the field of fractions of
$\Lambda$ and $\otimes=\otimes_\Lambda$. For a unitary $Q$-isomorphism 
$\varphi\colon H_1\otimes Q\to H_2\otimes Q$,
we define its {\em restricted graph} $\Gamma^0_\varphi$ by
\[
\Gamma^0_\varphi=\Gamma_\varphi\cap(H_1\oplus H_2)=\{h\oplus 
\varphi(h)\vert h\in H_1, \varphi(h)\in H_2\}
\subset H_1 \oplus H_2.
\]
If $\varphi$ is induced by a unitary $\Lambda$-isomorphism $f\colon 
H_1\to H_2$, then clearly $\Gamma^0_\varphi=\Gamma_f$.

\begin{lem}\label{lemma:symp}
Given any unitary isomorphism $\varphi\colon H_1\otimes Q\to 
H_2\otimes Q$, the restricted graph $\Gamma^0_\varphi$
is a Lagrangian relation $H_1\Rightarrow H_2$.
\end{lem}
\begin{pf}
Denote by $\omega_1$ (resp. $\omega_2$, $\omega$) the skew-hermitian 
form on $H_1$ (resp. $H_2$, $(-H_1)\oplus H_2$),
and pick $h,h'\in H_1$ such that $\varphi(h),\varphi(h')\in H_2$. 
Then,
\[
\omega(h\oplus \varphi(h),h'\oplus 
\varphi(h'))=-\omega_1(h,h')+\omega_2(\varphi(h),\varphi(h'))
=0.
\]
Therefore, $\Gamma^0_\varphi$ is isotropic. To check that it is 
Lagrangian, consider an element $x=x_1\oplus x_2$
of $\Ann(\Gamma^0_\varphi)\subset(-H_1)\oplus H_2$. For all $h$ in 
$H_1$ such that $\varphi(h)\in H_2$, 
\begin{eqnarray*}
0&=&
\omega(x,h\oplus\varphi(h))=-\omega_1(x_1,h)+\omega_2(x_2,\varphi(h))\\
&=&
-\omega_2(\varphi(x_1),\varphi(h))+\omega_2(x_2,\varphi(h))=
\omega_2(x_2-\varphi(x_1),\varphi(h)).
\end{eqnarray*}
Since $\varphi$ is an isomorphism, we have 
$H_2\subset\overline{\{\varphi(h)\vert h\in H_1, \varphi(h)\in 
H_2\}}$.
Therefore, $\omega_2(x_2-\varphi(x_1),h_2)=0$ for all $h_2\in H_2$. 
Since $\omega_2$ is non-degenerate, it follows
that $x_2=\varphi(x_1)$ so $x=x_1\oplus \varphi(x_1)\in 
\Gamma^0_\varphi$ and the lemma is proved.
\qed\end{pf}

Therefore, Lagrangian relations can be understood as a generalization 
of unitary isomorphisms. More precisely,
let ${\mathbf U_\Lambda}$ be the category of Hermitian 
$\Lambda$-modules and unitary $\Lambda$-isomorphisms. Also,
let ${\mathbf U^0_\Lambda}$ be the category of Hermitian 
$\Lambda$-modules, where the morphisms between $H_1$ and
$H_2$ are the unitary $Q$-isomorphisms between $H_1\otimes Q$ and 
$H_2\otimes Q$.

\begin{thm}
The maps $f\mapsto f\otimes id_Q$, $f\mapsto \Gamma_f$ and 
$\varphi\mapsto \Gamma^0_\varphi$ define embeddings of
categories
${\mathbf U_\Lambda}\subset{\mathbf U^0_\Lambda}$, 
${\mathbf U_\Lambda}\subset{\mathbf Lagr_\Lambda}$ and
${\mathbf U^0_\Lambda}\subset{\mathbf \overline{Lagr}_\Lambda}$ which 
fit in the commutative diagram
\[
\begin{CD}
{\mathbf U_\Lambda}@>>>{\mathbf U^0_\Lambda}\\
@VVV @VVV\\
{\mathbf Lagr_\Lambda}@>{\mathscr J}>>{\mathbf \overline{Lagr}_\Lambda}.
\end{CD}
\]
\end{thm}
\begin{pf}
The first embedding is clear. For the second one, note that the graph 
$\Gamma_f$ of a unitary $\Lambda$-isomorphism
$f$ is equal to the restricted graph of the induced unitary 
$Q$-isomorphism $f\otimes id_Q$. By Lemma \ref{lemma:symp},
$\Gamma_f$ is a Lagrangian relation. If $\Gamma_1$ and $\Gamma_2$ are 
the graphs of unitary $\Lambda$-isomorphisms
$f_1$ and $f_2$, then $\Gamma_2\Gamma_1$ is clearly the graph of 
$f_2\circ f_1$. Finally, two $\Lambda$-isomorphisms
with the same graph are equal.

By Lemma \ref{lemma:symp}, $\Gamma^0_\varphi$ is a Lagrangian 
relation, and it follows from the definition that
$\overline{\Gamma^0_\varphi}=\Gamma^0_\varphi$. Also, note that 
$\Gamma^0_\varphi\otimes Q=\Gamma_\varphi$.
Therefore, given two unitary $Q$-isomorphisms $\varphi_1$ and 
$\varphi_2$,
\begin{eqnarray*}
\Gamma^0_{\varphi_2\circ\varphi_1}&=&\Gamma_{\varphi_2\circ\varphi_1}
\cap(H_1\oplus H_3)=
\Gamma_{\varphi_2}\Gamma_{\varphi_1}\cap(H_1\oplus H_3)\\
&=&(\Gamma^0_{\varphi_2}\otimes Q)(\Gamma^0_{\varphi_1}\otimes 
Q)\cap(H_1\oplus H_3)=
(\Gamma^0_{\varphi_2}\Gamma^0_{\varphi_1}\otimes Q)\cap(H_1\oplus 
H_3)\\
&=&\overline{\Gamma^0_{\varphi_2}\Gamma^0_{\varphi_1}}=
\Gamma^0_{\varphi_2}\circ \Gamma^0_{\varphi_1}.
\end{eqnarray*}
It is clear that a $Q$-isomorphism $\varphi$ is entirely determined by 
its restricted graph $\Gamma^0_\varphi$.
\qed\end{pf}

\section{The Lagrangian representation}

\subsection{The category of oriented tangles}

Let $D^2$ be the closed unit disk in $\R^2$. Given a positive integer 
$n$, denote by $x_i$ the point $((2i-n-1)/n,0)$
in $D^2$, for $i=1,\dots,n$. Let $\eps$ and $\eps'$ be sequences of 
$\pm 1$ of respective length $n$ and $n'$.
An {\em $(\eps,\eps')$-tangle} is the pair consisting of the cylinder 
$D^2\times [0,1]$ and its oriented piecewise
linear $1$-submanifold $\tau$ whose oriented boundary $\partial\tau$ 
is $\sum_{j=1}^{n'}\eps'_j(x'_j,1)-\sum_{i=1}^{n}\eps_i(x_i,0)$.
Note that for such a tangle to exist, we must have 
$\sum_i\eps_i=\sum_j\eps'_j$.

Two $(\eps,\eps')$-tangles $(D^2\times [0,1],\tau_1)$ and
$(D^2\times[0,1],\tau_2)$ are {\em isotopic} if there exists an
auto-homeomorphism $h$ of $D^2\times [0,1]$, keeping 
$D^2\times\{0,1\}$ fixed, such that
$h(\tau_1)=\tau_2$ and $h\vert_{\tau_1}\colon\tau_1\simeq\tau_2$ is 
orientation-preserving. We shall denote by
$T(\eps,\eps')$ the set of isotopy classes of $(\eps,\eps')$-tangles, 
and by $id_\eps$ the isotopy class of the trivial
$(\eps,\eps)$-tangle $(D^2,\{x_1,\dots,x_n\})\times [0,1]$. 

\begin{figure}[Htb]
   \begin{center}
     \epsfig{figure=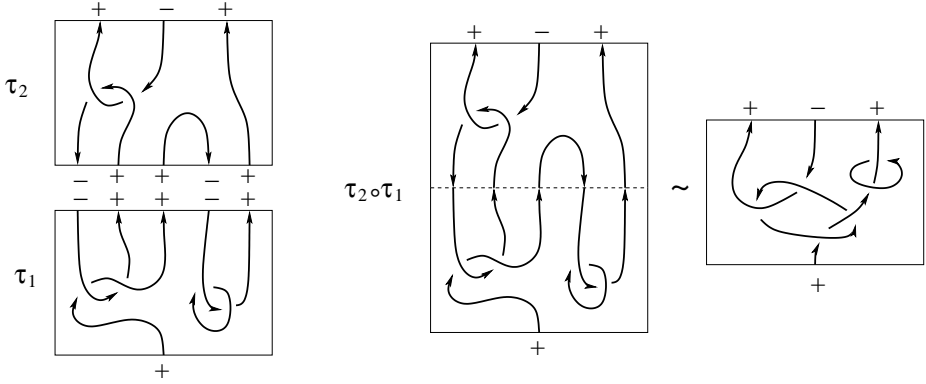,height=5cm}
     \caption{A tangle composition.}
     \label{fig:tangles}
   \end{center}
\end{figure}

Given an $(\eps,\eps')$-tangle $\tau_1$ and an $(\eps',\eps'')$-tangle 
$\tau_2$, their {\em composition} is the
$(\eps,\eps'')$-tangle $\tau_2\circ\tau_1$ obtained by gluing the two 
cylinders along the disk corresponding to $\eps'$
and shrinking the length of the resulting cylinder by a factor $2$ 
(see Figure \ref{fig:tangles}).
Clearly, the composition of tangles induces a composition
\[
T(\eps,\eps')\times T(\eps',\eps'')\longrightarrow T(\eps,\eps'')
\]
on the isotopy classes of tangles.

The {\em category of oriented tangles} ${\mathbf Tangles}$ is defined 
as follows: the objects are the finite sequences
$\eps$ of $\pm 1$, and the morphisms are given by 
$\hbox{Hom}(\eps,\eps')=T(\eps,\eps')$. The composition is clearly
associative, and the trivial tangle $id_\eps$ plays the role of the 
identity endomorphism of $\eps$. The aim of this
section is to construct a functor 
${\mathbf Tangles}\to{\mathbf Lagr_\Lambda}$.

\subsection{Objects}

\begin{figure}[Htb]
   \begin{center}
     \epsfig{figure=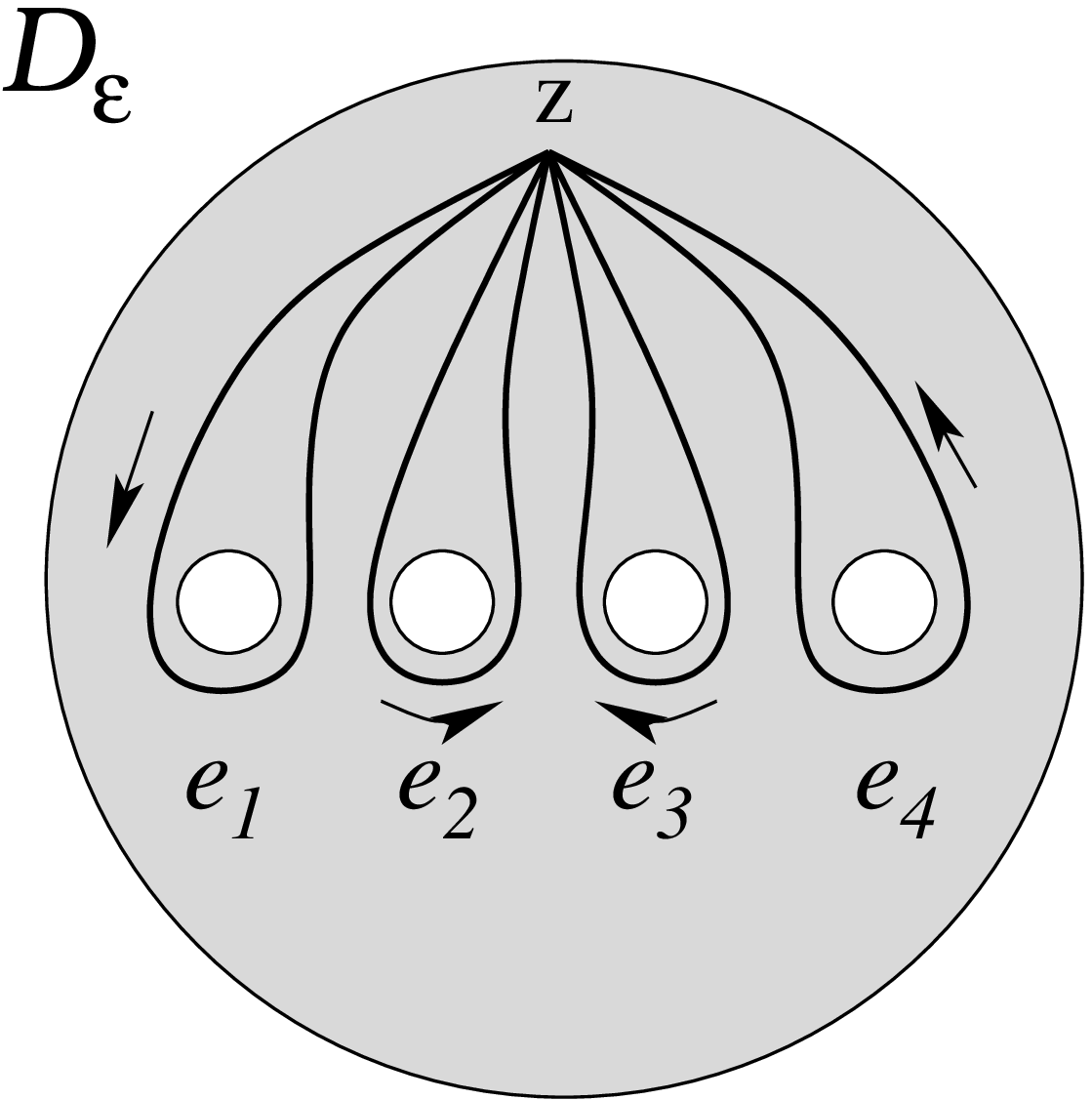,height=3.3cm}
     \caption{The space $D_\eps$ for $\eps=(+1,+1,-1,+1)$.}
     \label{fig:Deps}
   \end{center}
\end{figure}

Denote by $\N(\{x_1,\dots,x_n\})$ an open tubular neighborhood of 
$\{x_1,\dots,x_n\}$ in $D^2\subset\R^2$,
and by $S^2$ the $2$-sphere $\R^2\cup\{\infty\}$.
Given a sequence $\eps=(\eps_1,\dots,\eps_n)$ of $\pm 1$, let 
$\ell_\eps$ be the sum $\sum_{i=1}^n\eps_i$.
We shall denote by $D_\eps$ the compact surface
\[
D_\eps=\cases{D^2\setminus\N(\{x_1,\dots,x_n\})& if $\ell_\eps\neq 
0$;\cr
	S^2\setminus\N(\{x_1,\dots,x_n\})& if $\ell_\eps=0$,}
\]
endowed with the counterclockwise orientation, a base point $z$, and 
the generating family $\{e_1,\dots,e_n\}$
of $\pi_1(D_\eps,z)$, where $e_i$ is a simple loop turning once around 
$x_i$ counterclockwise if $\eps_i=+1$,
clockwise if $\eps_i=-1$ (see Figure \ref{fig:Deps}). The same space 
with the clockwise orientation will be denoted by $-D_\eps$.

The natural epimorphism $\pi_1(D_\eps)\to\Z$, $e_i\mapsto 1$ gives an 
infinite cyclic covering $\widehat{D}_\eps\to D_\eps$.
Choosing a generator $t$ of the group of the covering transformations 
endows the homology $H_1(\widehat{D}_\eps)$ with a
structure of module over $\Lambda=\Z[t,t^{-1}]$.
If $\ell_\eps\neq 0$, then $D_\eps$ retracts by deformation on the 
wedge of $n$ circles representing $e_1,\dots,e_n$,
and one easily checks that $H_1(\widehat{D}_\eps)$ is a free 
$\Lambda$-module with basis
$v_1=\hat e_1-\hat e_2,\dots,v_{n-1}=\hat e_{n-1}-\hat e_n$, where 
$\hat e_i$ is the path in $\widehat{D}_\eps$
lifting $e_i$ starting at some fixed lift $\hat z\in\widehat{D}_\eps$ 
of $z$.
If $\ell_\eps=0$, then $H_1(\widehat{D}_\eps)=
\bigoplus_i\Lambda v_i/\Lambda\hat\gamma$,
where $\hat\gamma$ is a lift of $\gamma=e_1^{\eps_1}\cdots 
e_n^{\eps_n}$ to $\widehat{D}_\eps$.
Note that in any case, $H_1(\widehat{D}_\eps)$ is a free 
$\Lambda$-module.

Let $<\;,\;>\colon H_1(\widehat{D}_\eps)\times 
H_1(\widehat{D}_\eps)\to\Z$ be the ($\Z$-bilinear, skew-symmetric)
intersection form induced by the orientation of $D_\eps$ lifted to 
$\widehat{D}_\eps$. Consider the pairing
$\omega_\eps\colon H_1(\widehat{D}_\eps)\times 
H_1(\widehat{D}_\eps)\to\Lambda$ given by
\[
\omega_\eps(x,y)=\sum_k<t^kx,y> t^{-k}.
\]
Note that this form is well-defined since, for any given $x,y\in 
H_1(\widehat{D}_\eps)$, the intersection $<t^kx,y>$
vanishes for all but a finite number of integers $k$. The 
multiplication by $t$ being an isometry with respect to the
intersection form, it is easy to check that $\omega_\eps$ is 
skew-hermitian with respect to the involution
$\Lambda\to\Lambda$ induced by $t\mapsto t^{-1}$.

\begin{exmp}
Consider $\eps$ of length $2$. If $\eps_1+\eps_2=0$, then 
$\widehat{D}_\eps$ is contractible so $H_1(\widehat{D}_\eps)=0$.
If $\eps_1+\eps_2\neq 0$, then $H_1(\widehat{D}_\eps)=\Lambda v$ with 
$v=\hat e_1-\hat e_2$, and
$\omega_\eps(v,v)=\frac{\eps_1+\eps_2}{2}(t-t^{-1})$, cf. Figure 
\ref{fig:disk}.
\begin{figure}[h]
   \begin{center}
     \epsfig{figure=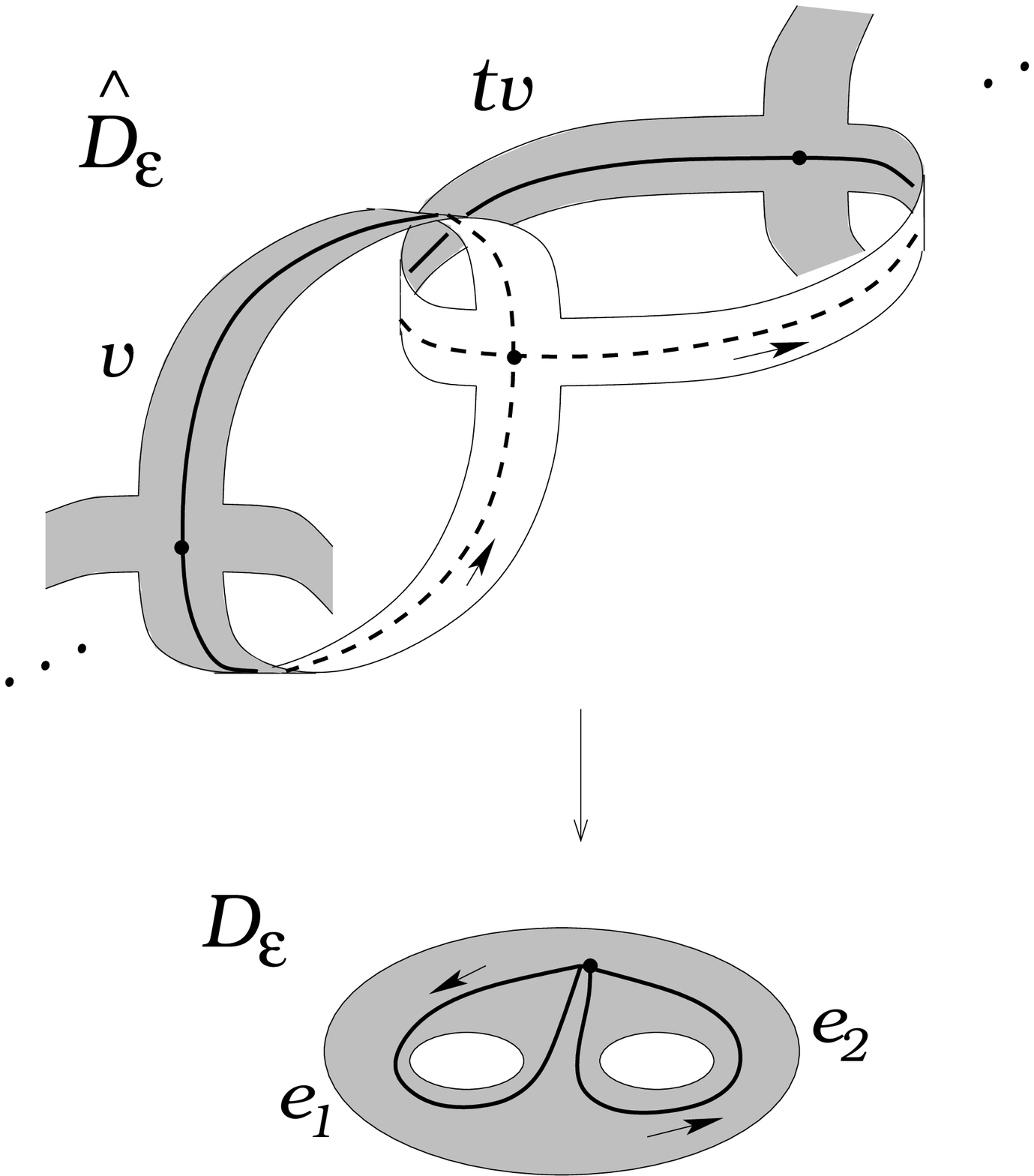,height=6cm}
     \caption{Computation of $\omega_\eps$ for $\eps=(+1,+1)$.}
     \label{fig:disk}
   \end{center}
\end{figure}
\end{exmp}
\noindent We shall give a proof of the following result in Section 4.

\begin{lem}\label{lemma:non_deg}
For any $\eps$, the form $\omega_\eps\colon 
H_1(\widehat{D}_\eps)\times H_1(\widehat{D}_\eps)\to\Lambda$ is 
non-degenerate.
\end{lem}

\subsection{Morphisms}

Given an $(\eps,\eps')$-tangle $\tau\subset D^2\times [0,1]$, denote 
by $\N(\tau)$ an open tubular neighborhood of $\tau$ and by $X_\tau$ 
its exterior
\[
X_\tau=\cases{(D^2\times [0,1])\setminus\N(\tau)& if
$\ell_\eps\neq 0$;\cr
	(S^2\times [0,1])\setminus\N(\tau)& if $\ell_\eps=0$.} 
\]
Note that $\ell_\eps=\ell_{\eps'}$.
We shall orient $X_\tau$ so that the induced orientation on $\partial 
X_\tau$ extends the orientation on $(-D_\eps)\sqcup D_{\eps'}$.
If $\ell_\eps\neq 0$, then the exact sequence of the pair
$(D^2\times [0,1],X_\tau)$ and the excision isomorphism give
\begin{eqnarray*}
H_1(X_\tau)&=&H_2(D^2\times 
[0,1],X_\tau)=H_2(\overline{\N(\tau)},\overline{\N(\tau)}\cap 
X_\tau)\\
	&=&\bigoplus_{j=1}^\mu 
H_2(\overline{\N(\tau_j)},\overline{\N(\tau_j)}\cap X_\tau),
\end{eqnarray*}
where $\tau_1,\dots,\tau_\mu$ are the connected components of $\tau$.
Since $(\overline{\N(\tau_j)},\overline{\N(\tau_j)}\cap X_\tau)$ is 
homeomorphic to $(\tau_j\times D^2,\tau_j\times S^1)$, we have 
$H_2(\overline{\N(\tau_j)},\overline{\N(\tau_j)}\cap X_\tau)=\Z m_j$, 
where $m_j$ is a meridian of $\tau_j$ oriented so that its linking 
number with $\tau_j$ is $1$.
Hence, $H_1(X_\tau)=\bigoplus_{j=1}^\mu \Z m_j$. If $\ell_\eps=0$,
then $H_1(X_\tau)=\bigoplus_{j=1}^\mu \Z m_j/\sum_{i=1}^n\eps_i e_i$.

The composition of the Hurewicz homomorphism and the homomorphism 
$H_1(X_\tau)\to\Z, m_j\mapsto 1$ gives an epimorphism 
$\pi_1(X_\tau)\to\Z$ which extends the previously defined 
homomorphisms $\pi_1(D_\eps)\to\Z$ and $\pi_1(D_{\eps'})\to\Z$.
As before, it determines an infinite cyclic covering 
$\widehat{X}_\tau\to X_\tau$, so the homology of $\widehat{X}_\tau$ is 
endowed with a natural structure of module over 
$\Lambda=\Z[t,t^{-1}]$.

Let $i_\tau\colon H_1(\widehat{D}_\eps)\to H_1(\widehat{X}_\tau)$ and 
$i'_\tau\colon H_1(\widehat{D}_{\eps'})\to H_1(\widehat{X}_\tau)$
be the homomorphisms induced by the obvious inclusion 
$\widehat{D}_\eps\sqcup\widehat{D}_{\eps'}\subset\widehat{X}_\tau$.
Denote by $j_\tau$ the homomorphism $H_1(\widehat{D}_\eps)\oplus 
H_1(\widehat{D}_{\eps'})\to H_1(\widehat{X}_\tau)$ given by 
$j_\tau(x,x')=i'_\tau(x')-i_\tau(x)$. Finally, set
\[
N(\tau)=\ker(j_\tau)\subset H_1(\widehat{D}_\eps)\oplus 
H_1(\widehat{D}_{\eps'}).
\]
Note that if $\tau$ and $\tau'$ are two isotopic 
$(\eps,\eps')$-tangles, then $N(\tau)=N(\tau')$.
 
\begin{lem}\label{lemma:Lag'}
$N(\tau)$ is a Lagrangian submodule of $(-H_1(\widehat{D}_\eps))\oplus 
H_1(\widehat{D}_{\eps'})$.
\end{lem}

\begin{lem}\label{lemma:comp}
If $\tau_1\in T(\eps,\eps')$ and $\tau_2\in T(\eps',\eps'')$, then 
$N(\tau_2\circ\tau_1)=N(\tau_2)N(\tau_1)$.
\end{lem}

We postpone the proof of these lemmas to the next section, and 
summarize our results in the following theorem.

\begin{thm}\label{thm:1}
Given a sequence $\eps$ of $\pm 1$, denote by ${\mathscr F}(\eps)$ the 
Hermitian $\Lambda$-module $(H_1(\widehat{D}_\eps),\omega_\eps)$.
For $\tau\in T(\eps,\eps')$, let ${\mathscr F}(\tau)$ be the 
Lagrangian relation
$N(\tau)\colon H_1(\widehat{D}_\eps)\Rightarrow 
H_1(\widehat{D}_{\eps'})$. Then, ${\mathscr F}$ is a functor
${\mathbf Tangles}\to{\mathbf Lagr_\Lambda}$.
\end{thm}

The usual notions of cobordism and $I$-equivalence for links 
generalize to
tangles in the obvious way. (The surface in $D^2\times [0,1]\times 
[0,1]$
interpolating between two tangles $\tau_1,\tau_2\subset D^2\times 
[0,1]$ should
be standard on $D^2\times\{0,1\}\times [0,1]$ and homeomorphic to 
$\tau_1\times [0,1]$.) It is easy to see (cf. \cite{KLW} Theorem 5.1 
and the proof of Proposition 5.3) that the Lagrangian relation 
$\overline{N(\tau)}$ is an $I$-equivalence invariant of $\tau$.

Finally, note that the usual computation of the Alexander module of a 
link $L$ from a diagram of $L$ extends to our setting. This gives a 
computation of $H_1(\widehat{D}_\eps)\oplus 
H_1(\widehat{D}_{\eps'})\stackrel{j_\tau}{\to}H_1(\widehat{X}_\tau)$ 
(cf. \cite[Proposition 4.4]{KLW}). Hence, it is possible to compute 
$N(\tau)$ from a diagram of $\tau$.

\section{Proof of the lemmas}

The proof of Lemmas \ref{lemma:non_deg} and \ref{lemma:Lag'} rely on 
the {\em Blanchfield duality theorem}. We recall this fundamental 
result referring for a proof and further details to
\cite[Appendix E]{Kaw}.

Let $M$ be a piecewise linear compact connected oriented  
$m$-dimensional manifold possibly with boundary. Consider an 
epimorphism of $\pi_1(M)$
onto a finitely generated free abelian group $G$. It induces a 
$G$-covering $\widehat M\to M$, so the homology modules of
$\widehat M$ are modules over $\Lambda=\Z G$.
For any integer $q$, let $<\;,\;>\colon H_q(\widehat M)\times 
H_{m-q}(\widehat M,\partial\widehat M)\to\Z$ be the $\Z$-bilinear 
intersection form
induced by the orientation of $M$ lifted to $\widehat M$. The
{\em Blanchfield pairing} is the  form $S\colon H_q(\widehat M)\times 
H_{m-q}(\widehat M,\partial\widehat M)\to\Lambda$ given by
\[
S(x,y)=\sum_{g\in G}<gx,y>g^{-1}.
\]
Note that $S$ is $\Lambda$-sesquilinear with respect to the involution 
of $\Lambda$ given by
$\sum_{g\in G}n_gg\mapsto\sum_{g\in G}n_gg^{-1}$.
The form $S$ induces a $\Lambda$-sesquilinear form
\[
H_q(\widehat M)/Tors_\Lambda H_q(\widehat M)\,\times 
\,H_{m-q}(\widehat M,\partial\widehat M)/Tors_\Lambda H_{m-q}(\widehat 
M,\partial\widehat M)\,\longrightarrow\Lambda.
\]

\begin{thm}[Blanchfield]
The latter form is non-degenerate.
\end{thm}

Let us now prove the lemmas stated in the previous section.

\begin{pf*}{Proof of Lemma \ref{lemma:non_deg}.}
Consider the Blanchfield pairing
\[
S_\eps\colon H_1(\widehat D_\eps)\times H_1(\widehat 
D_\eps,\partial\widehat D_\eps)\longrightarrow\Lambda.
\]
It follows from the definitions that 
$\omega_\eps(x,y)=S_\eps(x,j_\eps(y))$, where $j_\eps\colon 
H_1(\widehat D_\eps)\to H_1(\widehat D_\eps,\partial\widehat D_\eps)$ 
is the inclusion homomorphism.
Note that $\partial\widehat D_\eps$ consists of a finite number of 
copies of $\R$, so $H_1(\partial\widehat D_\eps)=0$ and $j_\eps$ is 
injective. Pick $y\in H_1(\widehat D_\eps)$ and assume that for all 
$x\in H_1(\widehat D_\eps)$, $0=\omega_\eps(x,y)=S_\eps(x,j_\eps(y))$. 
By the Blanchfield duality theorem, $j_\eps(y)\in 
Tors_\Lambda(H_1(\widehat D_\eps,\partial\widehat D_\eps))$, so 
$0=\lambda j_\eps(y)=j_\eps(\lambda y)$ for some $\lambda\in\Lambda$, 
$\lambda\neq 0$. Since $j_\eps$ is injective, $\lambda y=0$. As 
$H_1(\widehat D_\eps)$ is torsion-free, $y=0$, so $\omega_\eps$ is 
non-degenerate.
\qed\end{pf*}

\begin{pf*}{Proof of Lemma \ref{lemma:Lag'}.}
Let $H_1(\widehat{D}_\eps)\oplus 
H_1(\widehat{D}_{\eps'})\stackrel{i}{\to}H_1(\partial\widehat{X}_\tau)$
be the inclusion homomorphism, and denote by
\[
H_2(\widehat{X}_\tau,\partial\widehat{X}_\tau)\stackrel{\partial}
{\longrightarrow}H_1(\partial\widehat{X}_\tau)\stackrel{j}
{\longrightarrow}H_1(\widehat{X}_\tau)
\]
the homomorphisms appearing in the exact sequence of the pair 
$(\widehat{X}_\tau,\partial\widehat{X}_\tau)$. Also, denote by 
$\omega$ the pairing $(-\omega_\eps)\oplus\omega_{\eps'}$ on 
$(-H_1(\widehat{D}_\eps))\oplus H_1(\widehat{D}_{\eps'})$ and by
\[
S_{\partial X}\colon H_1(\partial\widehat{X}_\tau)\times 
H_1(\partial\widehat{X}_\tau)\to\Lambda\,,\quad
S_{X}\colon H_1(\widehat{X}_\tau)\times 
H_2(\widehat{X}_\tau,\partial\widehat{X}_\tau)\to\Lambda
\]
the Blanchfield pairings. Clearly, $N(\tau)=((-1)id\oplus id')(L)$, 
where $L=\ker(j\circ i)$ and $id$ (resp. $id'$) is the identity 
endomorphism of $H_1(\widehat{D}_\eps)$ (resp. 
$H_1(\widehat{D}_{\eps'})$). Then, $\Ann(N(\tau))=((-1)id\oplus 
id')\Ann(L)$
and we just need to check that $L$ is Lagrangian.

First, we check that $K=\ker(j)=\im(\partial)$ satisfies 
$\Ann_{\partial X}(K)=\overline K$,
where $\Ann_{\partial X}$ denotes the annihilator with respect to the 
form $S_{\partial X}$. 
Observe that for any $x\in  H_1( \partial\widehat{X}_\tau)$ and $Y\in 
H_2(\widehat{X}_\tau,\partial\widehat{X}_\tau)$,
we have $S_{\partial X}(x,\partial(Y))=S_X(j(x),Y)$. Therefore 
\begin{eqnarray*}
\Ann_{\partial X}(K)&=&\{x\in H_1(\partial\widehat{X}_\tau)\mid 
S_{\partial X}(x, K)=0\}\\
	&=&\{x\in H_1( \partial\widehat{X}_\tau)\mid 
S_{X}(j(x),H_2(\widehat{X}_\tau,\partial\widehat{X}_\tau))=0\}.
\end{eqnarray*}
By the Blanchfield duality, the latter set is just 
$j^{-1}(Tors_\Lambda(H_1(\widehat{X}_\tau)))=\overline K$.
 
Clearly, $i(L)\subset K$. The exact sequence of the pair 
$(\partial\widehat{X}_\tau,\widehat{D}_\eps\sqcup\widehat{D}_{\eps'})$ 
gives
\[
H_1(\widehat{D}_\eps)\oplus 
H_1(\widehat{D}_{\eps'})\stackrel{i}{\longrightarrow} 
H_1(\partial\widehat{X}_\tau)\longrightarrow T,
\]
where $T$ is a torsion $\Lambda$-module. This implies that 
$K\subset\overline{i(L)}$ and therefore $\overline{i(L)}=\overline K$. 
Since the forms $\omega$ and $S_{\partial X}$ are compatible under 
$i$,
\begin{eqnarray*}
\Ann(L)&=&i^{-1}(\Ann_{\partial X}(i(L)))=i^{-1}(\Ann_{\partial 
X}(\overline{i(L)}))=i^{-1}(\Ann_{\partial X}(\overline{K}))\\
&=&i^{-1}(\Ann_{\partial X}(K))=i^{-1}(\overline{K})=\overline L,
\end{eqnarray*}
so $L$ is Lagrangian and the lemma is proved.
\qed\end{pf*}

\begin{pf*}{Proof of Lemma \ref{lemma:comp}.}
Denote by $\tau$ the composition $\tau_2\circ\tau_1$.
Since $X_\tau=X_{\tau_1}\cup X_{\tau_2}$ and $X_{\tau_1}\cap 
X_{\tau_2}=D_{\eps'}$, we get the following Mayer-Vietoris exact 
sequence of $\Lambda$-modules:
\[
H_1(\widehat{D}_{\eps'})\stackrel{\alpha}{\to}H_1(\widehat 
X_{\tau_1})\oplus H_1(\widehat X_{\tau_2})\stackrel{\beta}{\to}
H_1(\widehat X_{\tau})\to 
H_0(\widehat{D}_{\eps'})\stackrel{\alpha_0}{\to}H_0(\widehat 
X_{\tau_1})\oplus H_0(\widehat X_{\tau_2}).
\]
The homomorphism $\alpha_0$ is clearly injective, so $\beta$ is onto 
and we get a short exact sequence which fits in the following 
commutative diagram
\[
\begin{CD}
0 @>>> H_{\eps'} @>{i}>> H_\eps\oplus H_{\eps'}\oplus H_{\eps''} 
@>{\pi}>> H_\eps\oplus H_{\eps''} @>>> 0\phantom{\,,}\\
& & @VV{\alpha}V @VV{\varphi}V @VV{j_{\tau}}V \\
0 @>>> \ker(\beta) @>>> H_1(\widehat X_{\tau_1})\oplus
H_1(\widehat X_{\tau_2}) @>{\beta}>> H_1(\widehat X_{\tau}) @>>>0\,,
\end{CD}
\]
where $H_\bullet$ denotes $H_1(\widehat{D}_\bullet)$, $i$ is the 
natural inclusion, $\pi$ the canonical projection, and 
$\varphi(x,x',x'')=(j_{\tau_1}(x,x'),j_{\tau_2}(x',x''))$. Clearly,
\[
\pi(\ker(\varphi))=\{x\oplus x''\mid \varphi(x,x',x'')=0\hbox{ for 
some $x'\in H_{\eps'}$}\}=\ker(j_{\tau_2})\ker(j_{\tau_1}).
\]
Therefore, we just need to check that 
$\pi(\ker(\varphi))=\ker(j_{\tau})$, which is an easy diagram chasing 
exercise using the surjectivity of $\alpha\colon 
H_{\eps'}\to\ker(\beta)$.
\qed\end{pf*}

\section{Examples}

\subsection{Braids}

An $(\eps,\eps')$-tangle $\tau=\tau_1\cup\dots\cup\tau_n\subset 
D^2\times [0,1]$ is called an
{\em oriented braid} if every component $\tau_j$ of $\tau$ is strictly 
increasing or strictly decreasing with respect to the projection to 
$[0,1]$.
Note that for such an oriented braid to exist, we must have 
$\sharp\{i\mid\eps_i=1\}=\sharp\{j\mid\eps'_j=1\}$ and 
$\sharp\{i\mid\eps_i=-1\}=\sharp\{j\mid\eps'_j=-1\}$. The finite 
sequences of $\pm 1$, as objects, and the isotopy classes of oriented 
braids, as morphisms, form a subcategory ${\mathbf Braids}$ of the 
category of oriented tangles. We shall now investigate the restriction 
of the functor $\mathscr F$ to this subcategory.
 
Consider an oriented braid $\beta=\beta_1\cup\dots\cup\beta_n\subset 
D^2\times [0,1]$. Clearly, there exists an isotopy $H_\beta\colon 
D^2\times [0,1]\to D^2\times [0,1]$ with $H_\beta(x,t)=(x,t)$ for 
$(x,t)\in(D^2\times\{0\})\cup(\partial D^2\times [0,1])$, such that 
$t\mapsto H_\beta(x_i,t)$ is
a homeomorphism of $[0,1]$ onto the arc $\beta_i$ for $i=1,\dots,n$. 
Let $h_\beta\colon D_\eps\to D_{\eps'}$ be the homeomorphism given by 
$x\mapsto H_\beta(x,1)$, and by the identity on $S^2\setminus D^2$ if 
$\eps_1+\dots+\eps_n=0$. It is a standard result that the isotopy 
class $(rel\;\partial D^2)$ of $h_\beta$ only depends on the isotopy 
class of $\beta$. Consider the lift $\hat h_\beta\colon \widehat 
D_\eps\to\widehat D_{\eps'}$ of $h_\beta$ fixing $\partial 
\widehat{D}^2$ pointwise, and denote by $f_\beta$ the induced unitary 
isomorphism $(\hat h_\beta)_\ast\colon H_1(\widehat{D}_\eps)\to 
H_1(\widehat{D}_{\eps'})$.

The isotopy $H_\beta$ provides a deformation retraction of $X_\beta$ 
to $D_{\eps'}$: let us identify $H_1(\widehat{X}_\beta)$ and 
$H_1(\widehat{D}_{\eps'})$ via this deformation. Clearly, the 
homomorphism $j_\beta\colon H_1(\widehat{D}_\eps)\oplus 
H_1(\widehat{D}_{\eps'})\to H_1(\widehat{X}_\beta)$ is given by 
$j_\beta(x,y)=y-f_\beta(x)$. Therefore,
\[
N(\beta)=\ker(j_\beta)=\{x\oplus f_\beta(x)\mid x\in 
H_1(\widehat{D}_\eps)\}=\Gamma_{f_\beta},
\]
the graph of the unitary isomorphism $f_\beta$. We have proved:

\begin{prop}\label{prop:braids}
The restriction of $\mathscr F$ to the subcategory of oriented braids 
gives a functor
${\mathbf Braids}\to{\mathbf U_\Lambda}$.
\end{prop}

Consider an $(\eps,\eps')$-tangle 
$\tau=\tau_1\cup\dots\cup\tau_n\subset D^2\times [0,1]$ such that 
every component
$\tau_i$ of $\tau$ is strictly increasing with respect to the 
projection to $[0,1]$.
Here, $\eps=\eps'=(1,\dots,1)$. We will simply call $\tau$ a
{\em braid}, or an {\em $n$-strand braid}.
As usual, we will denote by $B_n$ the group of isotopy classes of 
$n$-strand braids,
and by $\sigma_1,\dots,\sigma_{n-1}$ its standard set of generators 
(see Figure \ref{fig:elementary}).
Recall that the {\em Burau representation} $B_n\to GL_n(\Lambda)$ maps 
the generator $\sigma_i$ to the matrix
\[
I_{i-1}\oplus\pmatrix{1-t&t\cr 1&0}\oplus I_{n-i-1},
\]
where $I_k$ denotes the identity $(k\times k)$-matrix. This 
representation is reducible: it splits into the direct sum of an
$(n-1)$-dimensional representation $\rho$ and the trivial 
one-dimensional representation (see e.g. \cite{Bir}). 
Using the Artin presentation of $B_n$, one easily checks that the map
$\sigma_i\mapsto\rho(\sigma_i)^T$, where $(\quad)^T$ denotes the 
transposition, also defines a representation
$\rho^T\colon B_n\to GL_{n-1}(\Lambda)$.

\begin{prop}\label{prop:Burau}
The restriction of the functor $\mathscr F$ to $B_n$ gives a linear
anti-representation $B_n\to GL_{n-1}(\Lambda)$ which is the dual of 
$\rho^T$. 
\end{prop}
\begin{pf}
Consider two braids $\alpha,\beta\in B_n$. By Proposition 
\ref{prop:braids}, $N(\alpha)$ (resp. $N(\beta)$, $N(\alpha\beta)$)
is the graph of a unitary automorphism $f_\alpha$ (resp. $f_\beta$, 
$f_{\alpha\beta}$) of $H_1(\widehat{D}_\eps)$.
Note that the product $\alpha\beta \in B_n$ represents the composition 
$\beta\circ\alpha$ in the category of tangles.
Clearly, $f_{\alpha\beta}=f_\beta\circ f_\alpha$. Therefore,
$\mathscr F$ restricted to $B_n$ is an anti-representation.
In order to check that it corresponds to the dual of $\rho^T$,
we just need to verify that these anti-representations coincide on the 
generators $\sigma_i$ of $B_n$.
\begin{figure}[Htb]
   \begin{center}
     \epsfig{figure=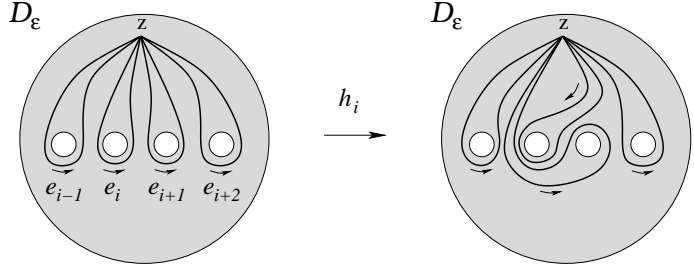,height=3.5cm}
     \caption{The action of $h_i$ on the loops 
$e_{i-1},\dots,e_{i+2}$.}
     \label{fig:Burau}
   \end{center}
\end{figure}

Denote by $f_i$ the unitary isomorphism corresponding to $\sigma_i$. 
We shall now compute the matrix of $f_i$ with respect to the basis 
$v_1,\dots,v_{n-1}$ of $H_1(\widehat{D}_\eps)$.
Consider the homeomorphism $h_i$ of $D_\eps$ associated with 
$\sigma_i$. As shown in Figure \ref{fig:Burau}, its action on the 
loops $e_j$ is given by
\[
h_i(e_j)=\cases{e_ie_{i+1}e_i^{-1} & if $j=i$;\cr
	e_i & if $j=i+1$;\cr
	e_j & else.}
\]
Therefore, the lift $\hat h_i$ of $h_i$ satisfies
\[
\hat h_i(\hat e_j)=\cases{\hat e_i-t(\hat e_i-\hat e_{i+1}) & if 
$j=i$;\cr
	\hat e_i & if $j=i+1$;\cr
	\hat e_j & else,}
\]
and the matrix of $f_i=(\hat h_i)_\ast$ with respect to the basis 
$v_j=\hat e_j-\hat e_{j+1}$ is
\begin{eqnarray*}
M_{f_1}&=&\pmatrix{-t&1\cr\phantom{-}0&1}\oplus I_{n-3},\quad 
M_{f_{n-1}}=I_{n-3}\oplus\pmatrix{1&\phantom{-}0\cr t&-t},\\
M_{f_i}&=&I_{i-2}\oplus
\pmatrix{1&0&0\cr t&-t&1\cr 0&0&1}\oplus 
I_{n-i-2}\;\;\hbox{ for }\;2\le i\le n-2.
\end{eqnarray*}
This is exactly $\rho(\sigma_i)$ (see, for instance, 
\cite[p.121]{Bir}).
\qed\end{pf}

\subsection{String links}

An $(\eps,\eps')$-tangle $\tau=\tau_1\cup\dots\cup\tau_n\subset 
D^2\times [0,1]$ is called an
{\em oriented string link} if every component $\tau_j$ of $\tau$ joins 
$D^2\times\{0\}$ and $D^2\times\{1\}$. Oriented string links clearly 
form a category ${\mathbf Strings}$ which satisfies
\[
{\mathbf Braids}\subset{\mathbf Strings}\subset{\mathbf Tangles},
\]
where all the inclusions denote embeddings of categories. Recall the 
functor ${\mathscr 
J}\colon{\mathbf Lagr_\Lambda}\to{\mathbf \overline{Lagr}_\Lambda}$
from Section 2.3.

\begin{prop}\label{prop:stringlinks}
The restriction of ${\mathscr J}\circ{\mathscr F}$ to the subcategory 
of oriented string links gives a functor
${\mathbf Strings}\to{\mathbf U^0_\Lambda}$.
\end{prop}
\begin{pf}
Since $\tau$ is an oriented string link, the inclusions $D_\eps\subset 
X_\tau$ and $D_{\eps'}\subset X_\tau$ induce isomorphisms in integral 
homology. Therefore, the induced homomorphisms $H_1(\widehat 
D_\eps;Q)\stackrel{i_\tau}{\to}H_1(\widehat X_\tau;Q)$ and 
$H_1(\widehat D_{\eps'};Q)\stackrel{i'_\tau}{\to}H_1(\widehat 
X_\tau;Q)$ are isomorphisms (see e.g. \cite[Proposition 2.3]{KLW}). 
Since $Q=Q(\Lambda)$ is a flat $\Lambda$-module,
$N(\tau)\otimes Q$ is the kernel of
\[
H_1(\widehat D_\eps;Q)\oplus H_1(\widehat 
D_{\eps'};Q)\stackrel{i'_\tau-i_\tau}{\longrightarrow}H_1(\widehat 
X_\tau;Q).
\]
Hence,
\[
{\mathscr J}\circ{\mathscr F}(\tau)=\overline{N(\tau)}=(N(\tau)\otimes 
Q)\cap(H_1(\widehat D_{\eps})\oplus H_1(\widehat 
D_{\eps'}))=\Gamma^0_\varphi,
\]
the restricted graph of the unitary $Q$-isomorphism 
$\varphi=(i'_\tau)^{-1}\circ i_\tau$.
\qed\end{pf}

If all the components of an oriented string link $\tau$ are oriented 
from bottom to top, we will simply speak of $\tau$ as a {\em string 
link}. By Proposition \ref{prop:stringlinks}, the restriction of
${\mathscr J}\circ{\mathscr F}$ to the category of string links gives 
a functor to the category ${\mathbf U^0_\Lambda}$. This functor is due 
to Le Dimet \cite{LeD} and was studied further in \cite{KLW}.

\subsection{Elementary tangles}

\begin{figure}[Htb]
   \begin{center}
     \epsfig{figure=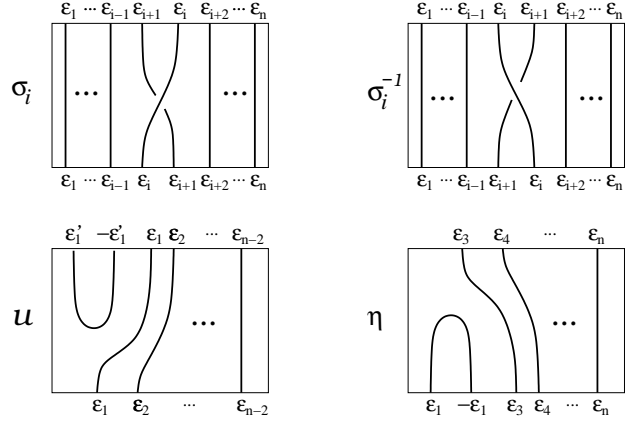,height=5.5cm}
     \caption{The elementary tangles.}
     \label{fig:elementary}
   \end{center}
\end{figure}

Every tangle $\tau\in T(\eps,\eps')$ can be expressed as a composition 
of the {\em elementary tangles} given in Figure \ref{fig:elementary}, 
where the orientation of the strands is determined by the signs $\eps$ 
and $\eps'$. We shall now compute explicity the functor
$\mathscr F$ on these tangles, assuming that $\ell_\eps\neq 0$.

Let us start with the tangle $u\in T(\eps,\eps')$. Here, 
$H_1(\widehat{D}_\eps)=\bigoplus_{i=1}^{n-3}\Lambda v_i$ and 
$H_1(\widehat{D}_{\eps'})=\bigoplus_{i=1}^{n-1}\Lambda v'_i$ where 
$v_i= \hat e_i-\hat e_{i+1}$ and $v'_i= \hat e'_i-\hat e'_{i+1}$. 
Moreover, $X_u$ is homeomorphic to the exterior of the trivial 
$(\eps'',\eps'')$-tangle, where 
$\eps''=(-\eps'_1,\eps_1,\dots,\eps_{n-2})=(\eps'_2,\dots,\eps'_n)$. 
Hence, $H_1(\widehat X_u)=\bigoplus_{i=1}^{n-2}\Lambda v''_i$ with 
$v''_i=\hat e''_i-\hat e''_{i+1}$ and the homomorphism $j_u\colon 
H_1(\widehat{D}_\eps)\oplus H_1(\widehat{D}_{\eps'})\to 
H_1(\widehat{X}_u)$ is given by $j_u(v_i)=-v''_{i+1}$ for 
$i=1,\dots,n-3$, $j_u(v'_1)=0$ and $j_u(v'_i)=v''_{i-1}$ for 
$i=2,\dots,n-1$. Therefore,
\[
N(u)=\ker(j_u)=\Lambda 
v'_1\oplus\bigoplus_{i=1}^{n-3}\Lambda(v_i\oplus v'_{i+2}).
\]
Similarly, we easily compute
\[
N(\eta)=\Lambda v_1\oplus\bigoplus_{i=1}^{n-3}\Lambda(v_{i+2}\oplus 
v'_i).
\]

Now, consider the oriented braid $\sigma_i\in T(\eps,\eps')$ given in 
Figure \ref{fig:elementary}. Then, $N(\sigma_i)$ is equal to the graph 
$\Gamma_{f_i}$ of a unitary isomorphism $f_i\colon 
H_1(\widehat{D}_\eps)\to H_1(\widehat{D}_{\eps'})$. As in the proof of 
Proposition \ref{prop:Burau}, we can compute the matrix $M_{f_i}$ of 
$f_i$ with respect to the bases $v_1,\dots,v_{n-1}$ of 
$H_1(\widehat{D}_\eps)$ and $v'_1,\dots,v'_{n-1}$ of 
$H_1(\widehat{D}_{\eps'})$:
\begin{eqnarray*}
M_{f_1}&=&\pmatrix{-t^{\eps_2}&1\cr 0&1}\oplus I_{n-3},\quad 
M_{f_{n-1}}=I_{n-3}\oplus\pmatrix{1&0\cr t^{\eps_n}&-t^{\eps_n}},\\
M_{f_i}&=&I_{i-2}\oplus\pmatrix{1&0&0\cr 
t^{\eps_{i+1}}&-t^{\eps_{i+1}}&1\cr 0&0&1}\oplus I_{n-i-2}\;\;
\hbox{ for }\;2\le i\le n-2.
\end{eqnarray*}
 
Finally, consider the tangle $\sigma^{-1}_i$ given in Figure 
\ref{fig:elementary}. Since it is an oriented braid, 
$N(\sigma^{-1}_i)$ is equal to the graph of a unitary isomorphism 
$g_i\colon H_1(\widehat{D}_{\eps'})\to H_1(\widehat{D}_\eps)$. 
Furthermore, we have
\[
\diag_{H_1(\widehat{D}_\eps)}=N(id_\eps)=N(\sigma^{-1}_i\circ\sigma_i)
=N(\sigma^{-1}_i)N(\sigma_i)=\Gamma_{g_i}\Gamma_{f_i}=
\Gamma_{g_i\circ f_i}.
\]
Therefore, $g_i\circ f_i$ is the identity endomorphism of 
$H_1(\widehat{D}_\eps)$, so the matrix of $g_i$ with respect to the 
basis given above is equal to $M_{g_i}=M^{-1}_{f_i}$.

With these elementary tangles, we can sketch an alternative
proof of 
Lemma \ref{lemma:Lag'} which does not make use of the Blanchfield 
duality. Indeed, any tangle $\tau\in T(\eps,\eps')$ can be written as 
a composition of $\sigma_i$, $\sigma^{-1}_i$, $u$ and $\eta$. By 
Lemmas \ref{lemma:Lag} and \ref{lemma:comp}, we just need to check 
that $N(\sigma_i)$, $N(\sigma^{-1}_i)$, $N(u)$ and $N(\eta)$ are 
Lagrangian. For $N(\sigma_i)$ and $N(\sigma^{-1}_i)$, this follows 
from Proposition \ref{prop:braids}, Lemma \ref{lemma:non_deg} and 
Lemma \ref{lemma:symp}. For $N(u)$ and $N(\eta)$, it can be verified 
by a direct computation of $\omega_\eps$.

\section{The module ${\mathbf N(\tau)}$}

\subsection{Freeness of ${\mathbf N(\tau)}$ and 
$\overline{{\mathbf N(\tau)}}$}

In this section, we deal with the following technical question: Given 
a tangle $\tau\in T(\eps,\eps')$, are the modules
$N(\tau)$ and $\overline{N(\tau)}$ free ? Clearly, these modules are 
contained in the free
module $H_1(\widehat{D}_\eps)\oplus H_1(\widehat{D}_{\eps'})$. But 
since the ring $\Z[t,t^{-1}]$ is not a principal
ideal domain, this is not sufficient to conclude that $N(\tau)$ and 
$\overline{N(\tau)}$ are free. Nevertheless, we have the following 
result.

Let us say that a tangle $\tau\in T(\eps,\eps')$ is {\em straight} if 
it has no closed components, and if at least one strand of $\tau$ 
joins $D_\eps$ with $D_{\eps'}$.

\begin{prop}\label{prop:free}
Given any tangle $\tau$, the $\Z[t,t^{-1}]$-module 
$\overline{N(\tau)}$ is free. If $\tau$ is a straight tangle, then 
$N(\tau)$ is also free.
\end{prop}

We shall need several notions of homological algebra, that we recall 
now. 
Let $\Lambda$ be a commutative ring with unit.
The {\em projective dimension pd(A)} of a $\Lambda$-module $A$ is the 
minimum integer $n$ (if it exists) such that there
is a projective resolution of length $n$ of $A$, that is, an exact 
sequence
\[
0\to P_n\to\cdots\to P_1\to P_0\to A\to 0
\]
where all the $P_i$'s are projective modules.
It is a well-known fact that if $0\to K_n\to P_{n-1}\to\cdots\to 
P_1\to P_0\to A\to 0$ is any resolution of $A$
with $pd(A)\le n$ and all the $P_i$'s projective, then $K_n$ is 
projective as well
(see, for instance, \cite[Lemma 4.1.6]{Wei}). The {\em global 
dimension} of a ring $\Lambda$ is the
(possibly $\infty$) number $\sup\{pd(A)\mid\hbox{$A$ is a 
$\Lambda$-module}\}$.
For example, the global dimension of $\Lambda$ is zero if $\Lambda$ is 
a field, and at most one if $\Lambda$
is a principal ideal domain.

Note that the ring $\Z[t,t^{-1}]$ has global dimension $2$ (see e.g. 
\cite[Theorem 4.3.7]{Wei}). We shall
also need the fact that all projective $\Z[t,t^{-1}]$-modules are free 
(\cite[Chapter 3.3]{Ros}). From now one, set
$\Lambda=\Z[t,t^{-1}]$.

\begin{lem}\label{lemma:tech}
Consider an exact sequence of $\Lambda$-modules
\[
0\longrightarrow K\longrightarrow P\longrightarrow F,
\]
where $P$ and $F$ are free $\Lambda$-modules. Then $K$ is free.
\end{lem}
\begin{pf}
Let $A$ be the image of the homomorphism $P\to F$. We claim that the 
projective dimension of $A$ is at most $1$.
Indeed, since the global dimension of $\Lambda$ is at most two, there 
is a projective resolution
$0\to P_2\stackrel{\partial}{\to}P_1\to P_0\to A\to 0$ of $A$. 
Splicing this resolution with the exact sequence
$0\to A\hookrightarrow F\to F/A\to 0$, we get a resolution of $F/A$
\[
0\to P_1/\partial P_2\to P_0\to F\to F/A\to 0,
\]
where $P_0$ and $F$ are projective. Since the global dimension of 
$\Lambda$ is $2$, we have $pd(F/A)\le 2$.
Hence, $P_1/\partial P_2$ is projective as well. Therefore, the 
resolution of $A$
\[
0\to P_1/\partial P_2\to P_0\to A\to 0
\]
is projective, so $pd(A)\le 1$. Now, the exact sequence $0\to K\to 
P\to A\to 0$ together with the fact that $P$ is free
and $pd(A)\le 1$, implies that $K$ is projective. Therefore, it is 
free.
\qed\end{pf}

\begin{lem}\label{lemma:free}
Let $H$, $H'$ and $H''$ be finitely generated free $\Lambda$-modules. 
Consider free submodules $N_1\subset H\oplus H'$
and $N_2\subset H'\oplus H''$ such that $(N_1\oplus N_2)\cap 
(0\oplus\diag_{H'}\oplus 0)=0$. Then $N_2N_1$ is a free submodule of
$H\oplus H''$.
\end{lem}
\begin{pf}
Denote by $f_1$ (resp.$f_1'$) the homomorphism $N_1\subset H\oplus 
H'\stackrel{\pi}{\to}H$
(resp. $N_1\subset H\oplus H'\stackrel{\pi'}{\to}H'$), where $\pi$ and 
$\pi'$ are the canonical projections.
Similarly, denote by $f'_2$ and $f''_2$ the homomorphisms $N_2\subset 
H'\oplus H''\to H'$ and $N_2\subset H'\oplus H''\to H''$.
Let $K$ be the kernel of $(-f'_1)\oplus f'_2\colon N_1\oplus N_2 \to 
H'$. 
  Our assumptions and Lemma \ref{lemma:tech} imply that $K$ is free. 
We have an exact sequence
\[
0\longrightarrow(N_1\oplus N_2)\cap (0\oplus\diag_{H'}\oplus 
0)\longrightarrow K
\stackrel{f_1\oplus f''_2}{\longrightarrow}N_2N_1\longrightarrow 0.
\]
Therefore, if $(N_1\oplus N_2)\cap (0\oplus\diag_{H'}\oplus 0)=0$, 
then 
$N_2N_1=K$ is free.
\qed\end{pf}

\begin{lem}\label{lemma:free'}
Consider tangles $\tau_1\in T(\eps,\eps')$, $\tau_2\in 
T(\eps',\eps'')$ such that $\tau_2\circ\tau_1$ is straight. Then
\[
(N(\tau_1)\oplus 
N(\tau_2))\cap(0\oplus\diag_{H_1(\widehat{D}_{\eps'})}\oplus 0)=0.
\]
\end{lem}
\begin{pf}
Denote by $\tau$ the tangle $\tau_2\circ\tau_1$. We claim that 
$H_2(X_\tau)=0$.
Let us first assume that $\ell_\eps\neq 0$. By excision,
\[
H_2(X_\tau)=H_3(D^2\times [0,1],X_\tau)=H_3(\tau\times D^2,\tau\times 
S^1)=0
\]
since $\tau$ has no closed components.
If $\ell_\eps=0$, consider the Mayer-Vietoris exact sequence 
associated with the decomposition
$X_\tau=((D^2\times[0,1])\setminus\N(\tau))\cup(D^2\times[0,1])$:
\[
0\to 
H_2(X_\tau)\to\Z\gamma\stackrel{i}{\to}H_1((D^2\times[0,1])\setminus\N
(\tau)),
\]
where $\gamma$ is a $1$-cycle parametrizing $\partial D^2$. Since one 
strand of $\tau$ joins $D_\eps$ with $D_{\eps''}$, we have 
$i(\gamma)\neq 0\in H_1((D^2\times[0,1])\setminus\N(\tau))=\Z^\mu$, 
where $\mu$ is the number of components of $\tau$.
Therefore, $i$ is injective, so $H_2(X_\tau)=0$ and the claim is 
proved.

Since $X_\tau$ has the homotopy type of a $2$-dimensional $CW$-complex 
and $H_2(X_\tau)=0$, we have $H_2(\widehat X_\tau)=0$. The 
decomposition $X_\tau=X_{\tau_1}\cup X_{\tau_2}$ gives the 
Mayer-Vietoris exact sequence
\[
H_2(\widehat X_\tau)=0\longrightarrow 
H_1(\widehat{D}_{\eps'})\stackrel{j}{\longrightarrow}H_1(\widehat 
X_{\tau_1})\oplus H_1(\widehat X_{\tau_2}).
\]
Therefore,
\begin{eqnarray*}
0&=&\ker(j)=\{x\in H_1(\widehat{D}_{\eps'})\mid j_{\tau_1}(0\oplus 
x)=j_{\tau_2}(x\oplus 0)=0\}\\
&\cong&(\ker(j_{\tau_1})\oplus\ker(j_{\tau_2}))
\cap(0\oplus\diag_{H_1(\widehat{D}_{\eps'})}\oplus 0)
\end{eqnarray*}
and the lemma is proved.
\qed\end{pf}

\begin{lem}\label{lemma:elementary}
Let $\tau$ be an elementary tangle, as described in Figure 
\ref{fig:elementary}.
Then, $N(\tau)$ is a free $\Lambda$-module.
\end{lem}
\begin{pf}
We already checked this statement by a direct computation when 
$\ell_\eps\neq 0$. In any case, note that $X_\tau$ has the homotopy 
type of a $1$-dimensional connected $CW$-complex $Y_\tau$ (unless 
$\tau$ is one of the $1$-strand tangles $u$ and $\eta$, in which case 
the lemma is obvious). Therefore, $H_1(\widehat X_\tau)$ is the kernel 
of ${\partial}\colon C_1(\widehat Y_\tau)\to C_0(\widehat Y_\tau)$. 
Since the latter two modules are free,
Lemma \ref{lemma:tech} implies that $H_1(\widehat X_\tau)$ is free. 
Now, consider the exact sequence
\[
0\to N(\tau)\hookrightarrow H_1(\widehat{D}_{\eps})\oplus 
H_1(\widehat{D}_{\eps'})\stackrel{j_\tau}{\to}H_1(\widehat X_\tau).
\]
Since $H_1(\widehat{D}_{\eps})\oplus H_1(\widehat{D}_{\eps'})$ and 
$H_1(\widehat X_\tau)$ are free, the conclusion follows from Lemma 
\ref{lemma:tech}.
\qed\end{pf}

\begin{pf*}{Proof of Proposition \ref{prop:free}.}
Consider the exact sequence
\[
0\to \overline{N(\tau)}\hookrightarrow H_1(\widehat{D}_{\eps})\oplus 
H_1(\widehat{D}_{\eps'})\to(H_1(\widehat{D}_{\eps})\oplus 
H_1(\widehat{D}_{\eps'}))/\overline{N(\tau)}\to 0.
\]
Clearly, the latter module is finitely generated and torsion free. 
Since $\Lambda$ is a noetherian ring, such a module embeds in a free 
$\Lambda$-module $F$, giving an exact sequence
\[
0\to \overline{N(\tau)}\hookrightarrow H_1(\widehat{D}_{\eps})\oplus 
H_1(\widehat{D}_{\eps'})\to F.
\]
By Lemma \ref{lemma:tech}, $\overline{N(\tau)}$ is free. The second 
statement follows from Lemmas \ref{lemma:comp}, \ref{lemma:free}, 
\ref{lemma:free'} and \ref{lemma:elementary}.
\qed\end{pf*}

Recall that for a $\Lambda$-module $N$, its rank $\rk_\Lambda N$ is 
defined by $\rk_\Lambda N=\dim_Q(N\otimes_\Lambda Q)$.
\begin{prop} \label{prop:rank}
Consider $\tau\in T(\eps,\eps')$ with $\eps$ of length $n$ and $\eps'$ 
of length $n'$. Then, the rank of $N(\tau)$ is given by
\[
\rk_\Lambda N(\tau)\;=\;\cases{0&if $n=n'=0$;\cr
				\frac{n+n'}{2}-1&if $\ell_\eps\neq 0$ 
or $nn'=0$ and $(n,n')\neq (0,0)$;\cr
				\frac{n+n'}{2}-2&if $\ell_\eps=0$ and 
$nn'>0$.}
\]
\end{prop}
\begin{pf}
Since $N(\tau)$ is a Lagrangian submodule of 
$H_1(\widehat{D}_\eps)\oplus H_1(\widehat{D}_{\eps'})$, we have
$\rk_\Lambda 
N(\tau)=\frac{1}{2}\rk_\Lambda(H_1(\widehat{D}_\eps)\oplus 
H_1(\widehat{D}_{\eps'}))$. If $\eps$ has length $n$, we know from 
Section 3.2
that
\[
\rk_\Lambda H_1(\widehat{D}_\eps)\;=\;\cases{0&if $n=0$;\cr
				n-1&if $\ell_\eps\neq 0$;\cr
				n-2&if $\ell_\eps=0$ and $n>0$.}
\]
The result follows.
\qed\end{pf}

\subsection{Recursive computation of ${\mathbf N(\tau)}$}

Consider two finitely generated free $\Lambda$-modules $H$ and $H'$ 
with fixed basis.
A homomorphism of $\Lambda$-modules $f\colon H\to H'$ is canonically 
described by its matrix $M_f$,
and the composition of homomorphisms corresponds to the product of 
matrices.
What about morphisms in the Lagrangian category? A free submodule $N$ 
of $H\oplus H'$ is determined by a matrix
of the inclusion $N\subset H\oplus H'$ with respect to a basis of $N$.
We will say that $N\subset H\oplus H'$ is {\em encoded} by this 
matrix.
For example, the graph $\Gamma_f$ of an isomorphism $f\colon H\to H'$ 
is encoded by the matrix $I\choose{M_f}$,
where $I$ is the identity matrix.

Let $H$, $H'$, $H''$ be finitely generated free $\Lambda$-modules with 
fixed basis.
Consider free submodules $N_1\subset H\oplus H'$ and $N_2\subset 
H'\oplus H''$.
A choice of a basis for $N_1$ and $N_2$ determines matrices 
${M_1}\choose{M_1'}$ and ${M'_2}\choose{M''_2}$ of the inclusions
$N_1\subset H\oplus H'$ and $N_2\subset H'\oplus H''$.
By Lemma \ref{lemma:free}, if
$(N_1\oplus N_2)\cap (0\oplus\diag_{H'}\oplus 0)=0$, then $N_2N_1$ is 
free. 
A natural question is: how can we compute a matrix
of the inclusion $N_2N_1\subset H\oplus H'$ from the matrices 
${M_1}\choose{M_1'}$ and ${M'_2}\choose{M''_2}$ ?

\begin{lem}\label{lemma:rec}
If $(N_1\oplus N_2)\cap (0\oplus\diag_{H'}\oplus 0)=0$, then the 
inclusion of $N_2N_1$ in $H\oplus H''$
is encoded by the matrix ${M_1W_1}\choose{M''_2W_2}$, where 
${W_1}\choose{W_2}$ is a matrix of the
inclusion of $K=\{x\in N_1\oplus N_2\mid (-M_1',M'_2)\cdot x=0\}$ in 
$N_1\oplus N_2$.
\end{lem}
\begin{pf}
We will assume the notation of the proof of Lemma \ref{lemma:free}.
By definition, $M_1$, $M_1'$, $M_2'$ and $M''_2$ are the matrices of 
$f_1$, $f'_1$, $f'_2$ and $f''_2$ with respect
to the bases of $N_1$, $N_2$, $H$, $H'$ and $H''$. Furthermore, we saw 
in the proof of Lemma \ref{lemma:free} that
$K=\ker((-f'_1)\oplus f'_2)$ is free. Let ${W_1}\choose{W_2}$ be a 
matrix of the inclusion
$K\subset N_1\oplus N_2$ with respect to a basis of $K$ and the fixed 
basis of $N_1\oplus N_2$. By definition,
$N_2N_1=(f_1\oplus f''_2)(K)$. Clearly, $\ker(f_1\oplus f''_2)\cap 
K=(N_1\oplus N_2)\cap (0\oplus\diag_{H'}\oplus 0)$.
Since the latter module is assumed to be trivial, $f_1\oplus f''_2$ 
restricted to $K$ gives an isomorphism onto $N_2N_1$.
The lemma follows easily.
\qed\end{pf}

Lemma \ref{lemma:rec} gives the following recursive method for the 
computation of $N(\tau)$, where $\tau$ is a straight tangle.

\begin{prop}\label{prop:rec}
Let $\tau_1\in T(\eps,\eps')$ and $\tau_2\in T(\eps',\eps'')$
be tangles such that $\tau_2\circ\tau_1$ is straight.
Then, $N(\tau_1)$, $N(\tau_2)$ and $N(\tau_2\circ\tau_1)$ are free.
Furthermore, if the inclusions $N(\tau_1)\subset 
H_1(\widehat{D}_\eps)\oplus H_1(\widehat{D}_{\eps'})$ and
$N(\tau_2)\subset H_1(\widehat{D}_{\eps'})\oplus 
H_1(\widehat{D}_{\eps''})$ are encoded by matrices
${M_1}\choose{M_1'}$ and ${M_2'}\choose{M_2''}$, then
$N(\tau_2\circ\tau_1)\subset H_1(\widehat{D}_\eps)\oplus 
H_1(\widehat{D}_{\eps''})$
is encoded by the matrix ${M_1W_1}\choose{M_2''W_2}$, where 
${W_1}\choose{W_2}$ is a matrix of the inclusion
$\{x\in N(\tau_1)\oplus N(\tau_2)\mid (-M_1',M'_2)\cdot x=0\}\subset 
N(\tau_1)\oplus N(\tau_2)$.
\end{prop}

In order to compute $N(\tau)$ for any straight tangle $\tau$, we only 
need to understand how $N(\tau_2\circ\tau_1)$ is
obtained from $N(\tau_1)$ for any elementary tangle $\tau_2$. Assuming 
the notation of Section 5.2, we easily get the following result.

\begin{prop}\label{prop:rec'}
Let $\tau_1\in T(\eps,\eps')$ be a straight tangle with $\ell_\eps\neq 
0$. If the inclusion $N(\tau_1)\subset H_1(\widehat{D}_\eps)\oplus 
H_1(\widehat{D}_{\eps'})$ is encoded by ${M_1}\choose{M_1'}$ and if 
$\tau_2\circ\tau_1$ is straight, then $N(\tau_2\circ\tau_1)$ is 
encoded by ${M\phantom{''}}\choose{M''}$, with
\begin{itemize}
\item{$M=(0\;M_1)$ and $M''={1\choose 0}\oplus M_1'$ if $\tau_2=u$;}
\item{$M=M_1$ and $M''=M_{f_i}^\epsilon M_1'$ if 
$\tau_2=\sigma_i^\epsilon$, for $\epsilon=\pm 1$;}
\item{$M=M_1W$ and $M''=\widetilde{M}_1'W$ if $\tau_2=\eta$, where 
$\widetilde{M}_1'$ denotes $M_1'$ without its first two lines
$\ell_1,\ell_2$, and $W$ encodes the solutions of the equation 
$\ell_2\cdot x=0$.}
\end{itemize}
\end{prop}

\section{The Alexander polynomial}

\begin{figure}[Htb]
   \begin{center}
     \epsfig{figure=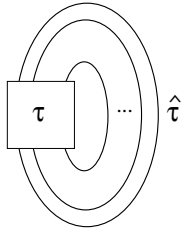,height=3cm}
     \caption{The closure $\hat\tau$ of an oriented tangle $\tau$.}
     \label{fig:closure}
   \end{center}
\end{figure}

Let $\tau\subset D^2\times [0,1]$ be an $(\eps,\eps)$-tangle, with 
$\eps$ of length $n$.
The {\em closure of $\tau$} is the oriented link $\hat\tau\subset S^3$ 
obtained from $\tau$
by adding $n$ oriented parallel strands in 
$S^3\setminus(D^2\times[0,1])$ as indicated in Figure 
\ref{fig:closure}.
The orientation of these strands is determined by $\eps$ in order to 
obtain a well-defined oriented link $\hat\tau$.

In this section, we show how the Alexander
polynomial $\Delta_{\hat\tau}$ of $\hat\tau$ is related to the 
Lagrangian module
$N(\tau)\subset H_1(\widehat{D}_\eps)\oplus H_1(\widehat{D}_\eps)$.

\subsection{Basics}

Let $\Lambda$ be a unique factorization domain. Consider a finite 
presentation $\Lambda^r\stackrel{f}{\to}\Lambda^g\to M\to 0$ of a 
$\Lambda$-module
$M$. We will denote by $\Delta(M)$ the greatest common divisor of the 
$(g\times g)$-minors of the matrix of $f$. It is well-known that, up 
to
multiplication by units of $\Lambda$, the element $\Delta(M)$ of 
$\Lambda$ only depends on the isomorphism class of $M$. Furthermore, 
if $0\to A\to B\to C\to 0$ is an exact sequence of $\Lambda$-modules, 
then $\Delta(B)\;\dot{=}\;\Delta(A)\Delta(C)$, where $\,\dot{=}\,$ 
denotes the
equality up to multiplication by units of $\Lambda$.

We briefly recall the definition of the $1$-variable Alexander 
polynomial of an oriented link $L\subset S^3$. Denote by $X_L$ the 
exterior of
$L$ in $S^3$, and consider the epimorphism $\pi_1(X_L)\to\Z$ given by 
the total linking number with $L$. It induces an infinite cyclic 
covering
$\widehat{X}_L\to X_L$. The $\Z[t,t^{-1}]$-module $H_1(\widehat{X}_L)$ 
is called the {\em Alexander module of $L$} and the Laurent polynomial 
$\Delta_L(t)=\Delta(H_1(\widehat{X}_L))$ is the {\em Alexander 
polynomial of $L$}. It is defined up to multiplication by $\pm t^\nu$, 
with $\nu\in\Z$.

\subsection{A factorization of the Alexander polynomial}

We use throughout this section the notation of Section 3.
\begin{lem}\label{lemma:Alex1}
For $\tau\in T(\eps,\eps)$ with $\ell_\eps\neq 0$,
\[
(t^{\ell_\eps}-1)\Delta_{\hat\tau}(t)\;\dot{=}\;(t-1)\Delta(A),
\]
where $A$ is the cokernel of $i'_\tau-i_\tau\colon 
H_1(\widehat{D}_\eps)\to H_1(\widehat{X}_\tau)$.
\end{lem}
\begin{pf}
Consider the compact manifold $Y_\tau$ obtained by pasting $X_\tau$ 
and $X_{id_\eps}$ along $D_\eps\sqcup D_\eps$. The epimorphisms
$\pi_1(X_\tau)\to\Z$ and $\pi_1(X_{id_\eps})\to\Z$ extend to an 
epimorphism $\pi_1(Y_\tau)\to\Z$ which defines a $\Z$-covering
$\widehat{Y}_\tau\to Y_\tau$. Hence, we have the Mayer-Vietoris exact 
sequence
\begin{eqnarray*}
H_1(\widehat{D}_\eps)\oplus 
H_1(\widehat{D}_\eps)&\stackrel{\alpha_1}{\to}&H_1(\widehat{D}_\eps)
\oplus H_1(\widehat{X}_\tau)
\stackrel{\beta}{\to}H_1(\widehat{Y}_\tau)\stackrel{\partial}{\to} 
H_0(\widehat{D}_\eps)\oplus H_0(\widehat{D}_\eps)\\
&\stackrel{\alpha_0}{\to}&H_0(\widehat{D}_\eps)\oplus 
H_0(\widehat{X}_\tau),
\end{eqnarray*}
where $\alpha_1(x,y)=(x+y,i_\tau(x)+i_\tau'(y))$. Since 
$H_0(\widehat{D}_\eps)=H_0(\widehat{X}_\tau)=\Lambda/(t-1)$,
the module $\im(\partial)=\ker(\alpha_0)$ is equal to $\Lambda/(t-1)$. 
This and the equality $A=\im(\beta)$ lead to the exact sequence
\[
0\to A\hookrightarrow H_1(\widehat{Y}_\tau)\to\Lambda/(t-1)\to 0,
\]
Hence, $\Delta(H_1(\widehat{Y}_\tau))\;\dot{=}\;(t-1)\Delta(A)$.

Clearly, $X_{\hat\tau}$ is the union of $Y_\tau$ and $D^2\times S^1$ 
along a torus $T\subset\partial Y_\tau$. The epimorphism 
$\pi_1(X_{\hat\tau}){\to}\Z$ given by the total linking number with 
$\hat\tau$ extends the previously defined epimorphism 
$\pi_1(Y_\tau)\to\Z$. Therefore, the exact sequence of the pair 
$(\widehat{X}_{\hat\tau},\widehat{Y}_\tau)$ gives
\[
0\to H_2(\widehat{Y}_\tau)\to 
H_2(\widehat{X}_{\hat\tau})\to\Lambda/(t^{\ell_\eps}-1)\to 
H_1(\widehat{Y}_\tau)\to H_1(\widehat{X}_{\hat\tau})\to 0.
\]
Note that both $H_2(\widehat{Y}_\tau)$ and 
$H_2(\widehat{X}_{\hat\tau})$ are free $\Lambda$-modules.
(This follows from the fact that $X_{\hat\tau}$ and $Y_\tau$ have the 
homotopy type of a $2$-dimensional $CW$-complex,
and from Lemma \ref{lemma:tech}.)
If $H_2(\widehat{X}_{\hat\tau})=0$, then
$\Delta(H_1(\widehat{Y}_\tau))\;\dot{=}\;(t^{\ell_\eps}-1)
\Delta(H_1(\widehat{X}_{\hat\tau}))=(t^{\ell_\eps}-1)
\Delta_{\hat\tau}(t)$ and the lemma holds.
If $H_2(\widehat{X}_{\hat\tau})\neq 0$, then 
$H_2(\widehat{Y}_\tau)\neq 0$ so both modules have positive rank. By 
an Euler characteristic
argument, the rank of $H_1(\widehat{X}_{\hat\tau})$ and 
$H_1(\widehat{Y}_\tau)$ is also positive. Therefore,
$\Delta(H_1(\widehat{X}_{\hat\tau}))=\Delta(H_1(\widehat{Y}_\tau))=0$, 
and the lemma is proved.
\qed\end{pf}

\begin{thm}\label{thm:Alex3}
Let $\tau\in T(\eps,\eps)$ be a tangle with $\ell_\eps\neq 0$, such 
that $N(\tau)$ is free. Then,
\[
\frac{t^{\ell_\eps}-1}{t-1}\,\Delta_{\hat\tau}(t)
\;\dot{=}\;\det(M'-M)\,\Delta(\cok(j_\tau)),
\]
where $M\phantom{'}\choose{M'}$ is a matrix of the inclusion 
$N(\tau)\subset H_1(\widehat{D}_\eps)\oplus H_1(\widehat{D}_\eps)$. 
\end{thm}
\begin{pf}
Since $N(\tau)=\ker(j_\tau)$, we have the exact sequence
\[
0\to N(\tau)\hookrightarrow H_1(\widehat{D}_\eps)\oplus 
H_1(\widehat{D}_\eps)\stackrel{j_\tau}{\to}H_1(\widehat{X}_{\tau})
\stackrel{\pi}{\to}\cok(j_\tau)\to 0.
\]
The module $A$ defined by the exact sequence 
$H_1(\widehat{D}_\eps)\stackrel{i'_\tau-i_\tau}{\longrightarrow}
H_1(\widehat{X}_\tau)\stackrel{p}{\to}A\to 0$
fits in the sequence
\[
N(\tau)\stackrel{\alpha}{\to}H_1(\widehat{D}_\eps)\stackrel{\beta}
{\to}A\stackrel{\gamma}{\to}\cok(j_\tau)\to 0,
\]
where $\alpha(x,y)=y-x$ for $x,y\in H_1(\widehat{D}_\eps)$, 
$\beta=p\circ i_\tau=p\circ i'_\tau$, and $\gamma(\zeta)=\pi(z)$ for 
$\zeta=p(z)\in A$,
$z\in H_1(\widehat{X}_\tau)$. We leave to the reader the proof that 
this sequence is exact. It then splits into two exact sequences
\[
N(\tau)\stackrel{\alpha}{\to}H_1(\widehat{D}_\eps)\stackrel{\beta}
{\to}\im(\beta)\to 0
\]
and
\[
0{\to}\im(\beta)\hookrightarrow A\to\cok(j_\tau)\to 0.
\]
The latter sequence implies that 
$\Delta(A)\;\dot{=}\;\Delta(\cok(j_\tau))\Delta(\im(\beta))$.
By Lemma \ref{lemma:Alex1}, we get
$(t^{\ell_\eps}-1)\Delta_{\hat\tau}(t)\;\dot{=}\;\Delta(\im(\beta))
\Delta(\cok(j_\tau))$. The
former sequence is nothing but a finite presentation of the module 
$\im(\beta)$. Furthermore, if a matrix of the inclusion 
$N(\tau)\subset H_1(\widehat{D}_\eps)\oplus H_1(\widehat{D}_\eps)$ is 
given by $M\phantom{'}\choose{M'}$, then a matrix of $\alpha$ is given 
by
$(M'-M)$. Since $N(\tau)$ is a Lagrangian submodule of 
$H_1(\widehat{D}_\eps)\oplus H_1(\widehat{D}_\eps)$, its rank is equal 
to the rank of
$H_1(\widehat{D}_\eps)$. Therefore, $M$ and $M'$ are square matrices 
and $\Delta(\im(\beta))\;\dot{=}\;\det(M'-M)$.
\qed\end{pf}

We have the following generalization of \cite[Theorem 3.11]{Bir}. 
(There, all the strands of the braid must be oriented in the same 
direction.)
\begin{cor}
If $\beta\in T(\eps,\eps)$ is an oriented braid with $\ell_\eps\neq 
0$, then
\[
\frac{t^{\ell_\eps}-1}{t-1}\,\Delta_{\hat\beta}(t)
\;\dot{=}\;\det(M_{f_\beta}-I),
\] 
where $M_{f_\beta}$ is a matrix of $f_\beta\colon 
H_1(\widehat{D}_\eps)\to H_1(\widehat{D}_\eps)$ (cf. Section 5.1) and 
$I$ is the identity matrix.
\end{cor} 
\begin{pf}
Since $N(\beta)=\Gamma_{f_\beta}$,
its inclusion in $H_1(\widehat{D}_\eps)\oplus H_1(\widehat{D}_\eps)$ 
is given by the matrix $I\choose{M_{f_\beta}}$.
Furthermore, $\widehat{D}_\eps$ is
a deformation retract of $\widehat{X}_\beta$, so the homomorphism 
$j_\beta$ is onto.
The equality then follows from Theorem \ref{thm:Alex3}.
\qed\end{pf}

A tangle $\tau\in T(\eps,\eps')$ is said to be {\em topologically 
trivial} if the oriented pair $(D^2\times [0,1],\tau)$ is homeomorphic 
to the
oriented pair $(D^2\times [0,1],id_{\eps''})$  for some $\eps''$. For 
instance, the oriented braids are topologically trivial, as well as 
the elementary tangles described in Figure \ref{fig:elementary}.
Note that a topologically tangle with $\ell_\eps\neq 0$ is always 
straight. Therefore, $N(\tau)$ is a free module if $\ell_\eps\neq 0$.

\begin{cor}\label{cor:toptrivial}
Consider a topologically trivial tangle $\tau\in T(\eps,\eps)$ with 
$\ell_\eps\neq 0$. Given a matrix
$M\phantom{'}\choose{M'}$ of the inclusion $N(\tau)\subset 
H_1(\widehat{D}_\eps)\oplus H_1(\widehat{D}_\eps)$, we have
\[
\delta\,\Delta_{\hat\tau}(t)\;\dot{=}\;\det(M'-M),
\]
where $\delta$ is a divisor of $\frac{t^{\ell_\eps}-1}{t-1}$ in 
$\Lambda$.
\end{cor}
\begin{pf}
Let $h$ be the homeomorphism between $(D^2\times[0,1],\tau)$ and 
$(D^2\times[0,1],id_{\eps''})$. The induced isomorphism
$h_\sharp\colon\pi_1(X_\tau)\to\pi_1(X_{id_{\eps''}})$ is compatible 
with the epimorphisms $\pi_1(X_\tau)\to\Z$ and 
$\pi_1(X_{id_{\eps''}})\to\Z$.
Therefore, $h$ lifts to a homeomorphism $\hat 
h\colon\widehat{X}_\tau\to\widehat{X}_{id_{\eps''}}$.

Denote by $B_\eps$ the compact surface $(\partial D^2\times[0,1])\cup 
(D_\eps\times\{0,1\})$. Since $\widehat{X}_{id_{\eps''}}$ retracts by 
deformation on
$\widehat{D}_{\eps''}\subset\widehat{B}_{\eps''}$, the manifold 
$\widehat{X}_\tau$ retracts by deformation on
$\widehat C=\hat h^{-1}(\widehat{D}_{\eps''})\subset\widehat{B}_\eps$. 
This leads to the following commutative diagram of inclusion 
homomorphisms,
\[
\xymatrix{
H_1(\widehat{D}_\eps)\oplus H_1(\widehat{D}_\eps) \ar[d]_i 
\ar[r]^{\quad\quad j\circ i}& H_1(\widehat{X}_\tau) \\ 
H_1(\widehat{B}_\eps)\ar[ur]^{j}
&\ar[l]_k H_1(\widehat{C})\ar[u]_{j\circ k} 
}
\]
where $j\circ k$ is an isomorphism. Let $\pi\colon 
H_1(\widehat{X}_\tau)\to\cok(j\circ i)$ and $\pi'\colon H_1(\widehat 
B_\eps)\to\cok(i)$ be the canonical
projections. Consider the homomorphism
$\varphi\colon\cok(j\circ i)\to\cok(i)$ given by
$\varphi(\pi(x))=\pi'\circ k\circ(j\circ k)^{-1}(x)$ for
$x\in H_1(\widehat X_\tau)$. We easily check that $\varphi$
is a well-defined injective homomorphism. Therefore, 
$\Delta(\cok(j_\tau))=\Delta(\cok(j\circ i))$ divides 
$\Delta(\cok(i))$. The exact sequence
of the pair 
$(\widehat{B}_\eps,\widehat{D}_\eps\sqcup\widehat{D}_\eps)$ gives
\[
H_1(\widehat{D}_\eps)\oplus 
H_1(\widehat{D}_\eps)\stackrel{i}{\to}
H_1(\widehat{B}_\eps)\to\Lambda/(t^{\ell_\eps}-1)
\to\Lambda/(t-1)\to 0.
\]
Therefore, $\Delta(\cok(i))\;\dot{=}\;(t^{\ell_\eps}-1)/(t-1)$. The 
result now follows from Theorem \ref{thm:Alex3}.
\qed\end{pf}

\subsection{Examples}

Given a topologically trivial tangle $\tau\in T(\eps,\eps)$ with 
$\ell_\eps\neq 0$, Propositions \ref{prop:rec}, \ref{prop:rec'} and 
Corollary \ref{cor:toptrivial} provide a method for the computation of 
the Alexander polynomial $\Delta_{\hat\tau}$.
We now give several examples of such computations.
\bigskip

\begin{figure}[Htb]
   \begin{center}
     \epsfig{figure=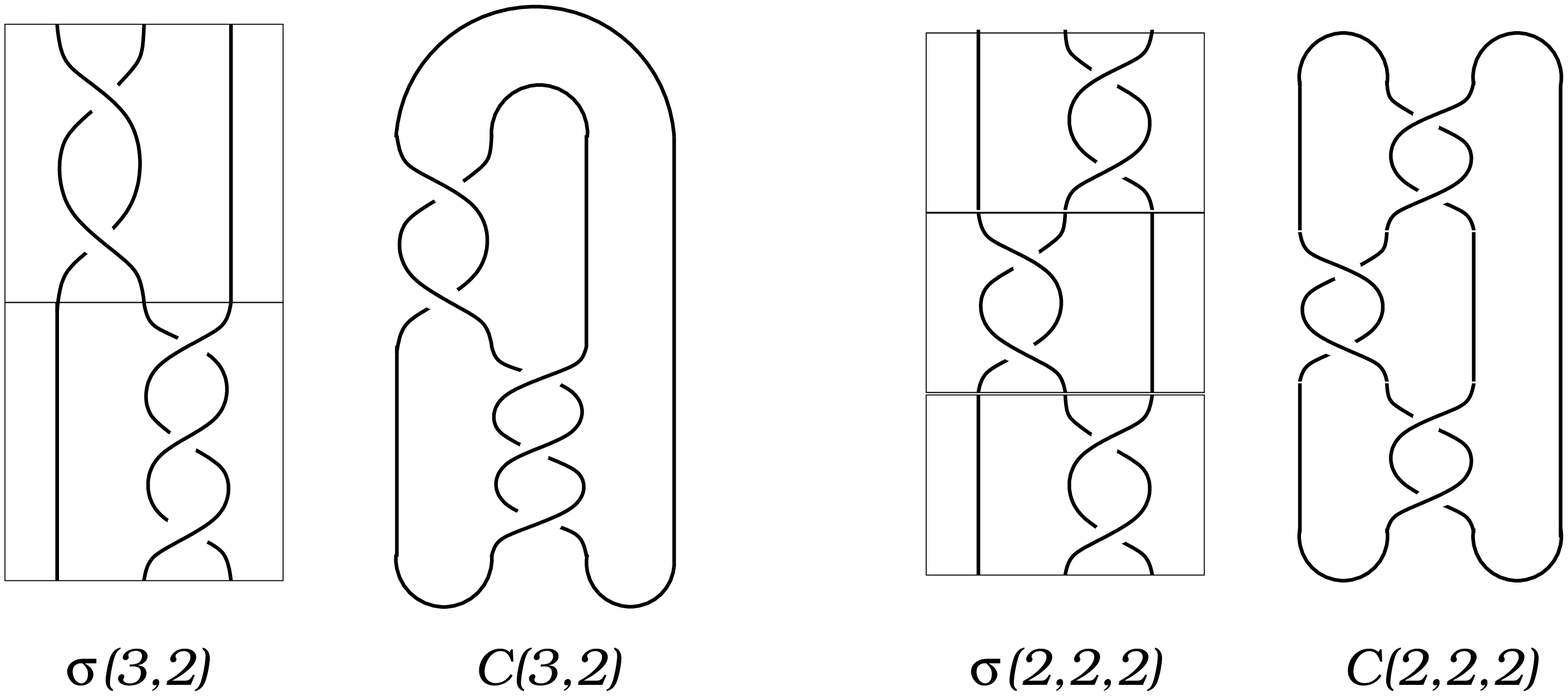,height=4cm}
     \caption{Rational tangles and rational links.}
     \label{fig:rational}
   \end{center}
\end{figure}

\noindent{\bf Rational links.} For integers $a_1,\dots,a_n$, denote by 
$\sigma(a_1,\dots,a_n)$ the following unoriented $3$-strand braid:
\[
\sigma(a_1,\dots,a_n)=\cases{\sigma_2^{a_1}\sigma_1^{-a_2}
\sigma_2^{a_3}\cdots\sigma_2^{a_n}& if $n$ is odd;\cr
\sigma_2^{a_1}\sigma_1^{-a_2}\sigma_2^{a_3}\cdots\sigma_1^{-a_n}& if 
$n$ is even.}
\]
Consider the unoriented $3$-strand tangle
$\tau(a_1,\dots,a_n)=\tau_n\circ \sigma(a_1,\dots,a_n)$, where
\[
\tau_n=\cases{u\circ\eta& if $n$ is odd;\cr 
u\circ\eta\circ\sigma_2\circ\sigma_1& if $n$ is even.}
\]
(Recall Figure \ref{fig:elementary} for the definition of the tangles 
$u$, $\eta$ and $\sigma_i$.)
Finally, denote by $C(a_1,\dots,a_n)$ the unoriented link given by the 
closure of $\tau(a_1,\dots,a_n)$. Such a link is 
called a {\em rational link} or a {\em $2$-bridge link} (see 
\cite{Con} and Figure \ref{fig:rational} for examples).

Consider the oriented link $L$ obtained by endowing $C(a_1,\dots,a_n)$ 
with an orientation. (Note that there is no canonical way
to do so: $L$ is not uniquely determined by the integers 
$(a_1,\dots,a_n)$.) This turns $\sigma(a_1,\dots,a_n)$ into
an oriented braid $\beta$, and one easily compute the associated 
matrix
$M_{f_\beta}=\pmatrix{m_{11}&m_{12}\cr m_{21}&m_{22}}$, where 
$m_{ij}\in\Lambda$ for $i,j=1,2$.
\begin{prop}\label{prop:rational}
The Alexander polynomial of $L$ is given by
\[
\Delta_L(t)\;\dot{=}\;\cases{m_{21}& if $n$ is odd;\cr m_{11}& if $n$ 
is even.}
\]
\end{prop}
\begin{pf}
Let us first assume that $n$ is odd. Consider the decomposition 
$\tau=\tau_n\circ\beta$. In the canonical bases $v_1,v_2$ of
$H_1(\widehat{D}_\eps)$ and $v'_1,v'_2$ of $H_1(\widehat{D}_{\eps'})$, 
the inclusion
$N(\tau_n)\subset H_1(\widehat{D}_{\eps'})\oplus 
H_1(\widehat{D}_\eps)$ is encoded by the matrix
${M'}\choose{M\phantom{'}}$, with $M'=\pmatrix{1&0\cr 0&0}$ and 
$M=\pmatrix{0&1\cr0&0}$. Furthermore,
the inclusion $N(\beta)\subset H_1(\widehat{D}_{\eps})\oplus 
H_1(\widehat{D}_{\eps'})$ is encoded by the matrix
${I}\choose{M_{f_\beta}}$. Since $M_{f_\beta}$ is invertible, the 
solutions of the system $(-M_{f_\beta},M')\cdot x=0$ are given by
${{W_1}\choose{W_2}}={{M^{-1}_{f_\beta}M'}\choose{I}}$. By Proposition 
\ref{prop:rec}, $N(\tau)$ is encoded by
${M^{-1}_{f_\beta}M'}\choose{M}$. By Corollary \ref{cor:toptrivial},
\begin{eqnarray*}
\Delta_L(t)&\dot{=}&\det(M-M^{-1}_{f_\beta}M')\;\dot{=}\;
\det(M_{f_\beta}M-M')\\
	&=&\det\left(\pmatrix{m_{11}&m_{12}\cr m_{21}&m_{22}} 
\pmatrix{0&1\cr0&0} - 	\pmatrix{1&0\cr0&0}\right)\;\dot{=}\;m_{21}.
\end{eqnarray*}
If $n$ is even, we have $M'=\pmatrix{0&0\cr 1&0}$ and 
$M=\pmatrix{0&1\cr0&0}$. This leads to $\Delta_L(t)\;\dot{=}\;m_{11}$.
\qed\end{pf}

For example, consider an oriented knot $K$ obtained by
orienting the knot $C(3,2)$ 
described in Figure \ref{fig:rational}. The corresponding oriented 
braid $\beta$ is the composition of $5$ elementary braids, leading to
\begin{eqnarray*}
M_{f_\beta}&=&\pmatrix{-t^\epsilon &1\cr 0&1}^{-2}\pmatrix{1 &0\cr 
t^\epsilon&-t^\epsilon}\pmatrix{1 &0\cr 
t^{-\epsilon}&-t^{-\epsilon}}\pmatrix{1 &0\cr 
t^\epsilon&-t^\epsilon}\\
&=&\pmatrix{2t^{-2\epsilon}-3t^{-\epsilon}+2&t^{-\epsilon}-1\cr 
2t^\epsilon-1&-t^\epsilon},
\end{eqnarray*}
where $\epsilon$ is $\pm 1$ according to the orientation of $K$. By 
Proposition \ref{prop:rational}, we have 
$\Delta_K(t)\;\dot{=}\;2t-3+2t^{-1}$.

Let $L$ be an oriented link obtained by orienting $C(2,2,2)$ so that 
the linking number of the components is $+2$. Here, we get
\begin{eqnarray*}
M_{f_\beta}&=&\pmatrix{1&0\cr t^\epsilon 
&-t^\epsilon}^2\pmatrix{-t^\epsilon&1\cr 
0&1}^{-1}\pmatrix{-t^{-\epsilon}&1\cr 0&1}^{-1}\pmatrix{1 &0\cr 
t^\epsilon&-t^\epsilon}^2\\
&=&\pmatrix{t^{2\epsilon}-2t^{\epsilon}+2&
t^{\epsilon}-t^{2\epsilon}\cr 
2(t^\epsilon-t^{2\epsilon})(t^{2\epsilon}-t^\epsilon+1)
&2t^{4\epsilon}
-2t^{3\epsilon}+t^{2\epsilon}},
\end{eqnarray*}
where $\epsilon=\pm 1$ depends on the global orientation of $L$. 
Therefore, $\Delta_L(t)\;\dot{=}\;2(t-1)(t-1+t^{-1})$.
Finally, if we orient $C(2,2,2)$ so that the linking number of the 
components is $-2$, the resulting oriented link $L'$ has Alexander 
polynomial $\Delta_{L'}\;\dot{=}\;(t-1)(t-4+t^{-1})$.

\noindent{\bf $\mathbf 2$-strand tangles.} In this section, we use the 
techniques introduced above to define an invariant
of $(2,2)$-tangles formed by two arcs and having no closed components.
This invariant is a pair of elements of $\Lambda$ defined up to 
simultaneous
multiplication by a unit of $\Lambda$. We study the behaviour of this 
invariant
under the basic transformations of $(2,2)$-tangles introduced by 
Conway \cite{Con}.

Consider a tangle $\tau\in T(\eps,\eps')$ with no closed components, 
where $\eps$ and $\eps'$ are sequences of length $2$.
By bending $\tau$, we get a tangle $\tau^b\in T(\emptyset,\mu)$ where 
$\emptyset$ is the empty sequence and 
$\mu=(\eps'_1,\eps'_2,-\eps_2,-\eps_1)$.

\begin{lem}\label{lemma:2-strand}
The submodule $N(\tau^b)$ of $H_1(\widehat D_\mu)$ is free of rank 
one.
\end{lem}
\begin{pf}
One can write $\tau^b$ as a composition $\tau^b=\tau'\circ u$, where 
$u\in T(\emptyset,\tilde\eps)$ is the elementary $1$-strand `cup' 
tangle and
$\tau'\in T(\tilde\eps,\eps)$ is a straight tangle.
Since $H_1(\widehat D_\emptyset)=H_1(\widehat D_{\tilde\eps})=0$, we 
have $N(u)=0$. By Lemma \ref{lemma:comp}, $N(\tau^b)=N(\tau')$ which 
is free
by Proposition \ref{prop:free}. Its rank is one by Proposition 
\ref{prop:rank}.
\qed\end{pf}
Recall from Section 3.2 that $H_1(\widehat D_\mu)=(\Lambda 
v_1\oplus\Lambda v_2\oplus\Lambda v_3)/\Lambda\hat\gamma$,
where $v_i=\hat e_i-\hat e_{i+1}$ and $\gamma=e_1^{\eps_1}\cdots 
e_4^{\eps_4}$. Therefore, $H_1(\widehat D_\mu)$ is free with basis 
$v_1,v_2$.
Using this fact and Lemma \ref{lemma:2-strand}, the inclusion 
$N(\tau^b)\subset H_1(\widehat D_\mu)$ is given by a matrix
$\pmatrix{m_1\cr m_2}$ with $m_1,m_2\in\Lambda$,
unique up to multiplication by $\pm t^\nu$ with $\nu\in\Z$. Let us 
denote this by $\tau\sim\pmatrix{m_1,m_2}$.

\begin{figure}[Htb]
   \begin{center}
     \epsfig{figure=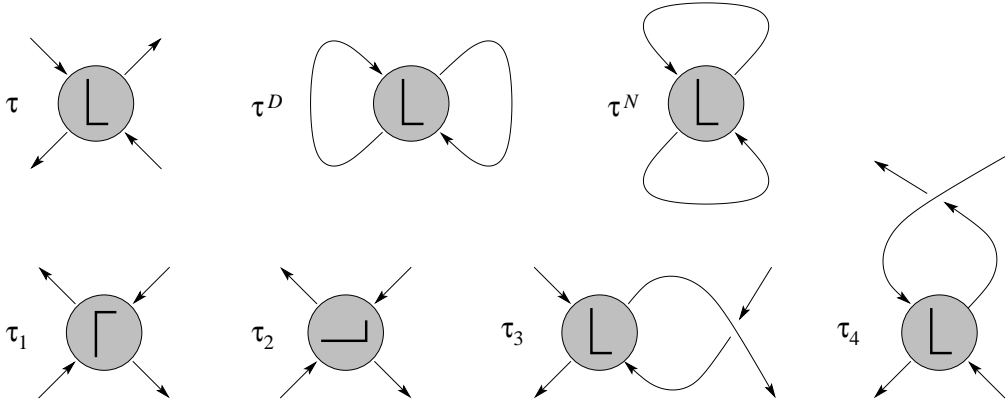,height=5.3cm}
     \caption{Tangles with two strands.}
     \label{fig:2-tangles}
   \end{center}
\end{figure}

For concreteness, we shall assume throughout the rest of the 
discussion that $\eps=\eps'=(-1,+1)$ as for the tangle $\tau$ in 
Figure \ref{fig:2-tangles}. (The other five cases can be treated 
similarly.)
Consider the tangles $\tau_1,\tau_2,\tau_3$ and $\tau_4$ shown in 
Figure \ref{fig:2-tangles}: $\tau_1$ is obtained from
$\tau$ by a horizontal reflexion, $\tau_2$ by a rotation to the angle 
$\pi/2$, $\tau_3$ by addition of a twist to the right, and
$\tau_4$ by addition of a twist to the top.

\begin{prop}\label{prop:move}
If $\tau\sim \pmatrix{m_1,m_2}$, then $\tau_1\sim\pmatrix{m_1,-m_2}$, 
$\tau_2\sim\pmatrix{m_2,-m_1}$, $\tau_3\sim\pmatrix{tm_1,m_1-m_2}$ and 
$\tau_4\sim\pmatrix{m_2-tm_1,m_2}$.
\end{prop}
\begin{pf}
We have $\tau^b\in T(\emptyset,\mu)$ with $\mu=(-1,+1,-1,+1)$, while 
$\tau_1^b,\tau_2^b\in T(\emptyset,\mu')$ where
$\mu'=(+1,-1,+1,-1)$. Hence, $H_1(\widehat D_{\mu'})=(\Lambda 
v'_1\oplus\Lambda v'_2\oplus\Lambda v'_3)/\Lambda(v'_1+v'_3)$.
The horizontal reflexion induces an isomorphism $H_1(\widehat 
D_\mu)\to H_1(\widehat D_{\mu'})$ given by $v_1\mapsto -v'_3=v'_1$ and 
$v_2\mapsto -v'_2$.
Hence, $\tau_1\sim\pmatrix{m'_1,m'_2}$ with $\pmatrix{m'_1\cr 
m'_2}=\pmatrix{1&\phantom{-}0\cr 0&-1}\pmatrix{m_1\cr 
m_2}=\pmatrix{\phantom{-}m_1\cr -m_2}$.
Similary, the rotation to the angle $\pi/2$ induces an isomorphism 
$H_1(\widehat D_\mu)\to H_1(\widehat D_{\mu'})$ given by the matrix 
$\pmatrix{\phantom{-}0&1\cr -1&0}$. Thus, 
$\tau_2\sim\pmatrix{m_2,-m_1}$.
Note that $\tau^b_3\in T(\emptyset,\mu'')$, where 
$\mu''=(-1,-1,+1,+1)$.
The transformation from $\tau^b$ to $\tau^b_3$ can be understood as a 
composition $\tau^b_3=\sigma\circ \tau^b$, where
$\sigma$ is a spherical braid.
By the results of Section 5.3, the isomorphism $H_1(\widehat D_\mu)\to 
H_1(\widehat D_{\mu''})$ corresponding to $\sigma$ is given by
$v_1\mapsto v''_1+t^{-1}v''_2$ and $v_2\mapsto -t^{-1}v''_2$. 
Therefore, $\tau_3\sim\pmatrix{m_1,t^{-1}(m_1-m_2)}$, which is 
equivalent to
$\pmatrix{tm_1,m_1-m_2}$. The case of $\tau_4$ is similar.
\qed\end{pf}

\begin{prop}
If $\tau$ is topologically trivial and $\tau\sim\pmatrix{m_1,m_2}$, 
then the oriented links $\tau^D$ and $\tau^N$ described in Figure 
\ref{fig:2-tangles}
have the Alexander module
\[
H_1(\widehat X_{\tau^D})=\Lambda/(m_1)\quad\hbox{ and }\quad 
H_1(\widehat X_{\tau^N})=\Lambda/(m_2).
\]
In particular, $\Delta_{\tau^D}(t)\;\dot{=}\;m_1$ and 
$\Delta_{\tau^N}(t)\;\dot{=}\;m_2$.
\end{prop}
\begin{pf}
Since $\tau$ is topologically trivial, $H_1(\widehat 
X_{\tau^b})=H_1(\widehat X_\tau)=\Lambda$ and the inclusion 
homomorphism
$j\colon H_1(\widehat D_\mu)=\Lambda v_1\oplus \Lambda v_2\to 
H_1(\widehat X_\tau)$ is onto (cf. the proof of Corollary 
\ref{cor:toptrivial}).
Therefore, the greatest common divisor of $j(v_1)$ and $j(v_2)$ is 
$1$. Hence, the kernel $N(\tau^b)$ of $j$ is generated by 
$j(v_2)v_1-j(v_1)v_2$, so $m_1=j(v_2)$ and $m_2=-j(v_1)$. Since the 
exterior of $\tau^D$ in $S^3$ can be written $X_{\tau^D}=X_\tau\cup 
X_{id}$, we have the Mayer-Vietoris exact sequence
\[
H_1(\widehat D_\mu)\stackrel{\varphi}{\to}H_1(\widehat X_\tau)\oplus 
H_1(\widehat X_{id})\to H_1(\widehat X_{\tau^D})\to 0.
\]
Clearly, $H_1(\widehat X_{id})=\Lambda v_1$ and a matrix of $\varphi$ 
is given by $\pmatrix{j(v_1)&j(v_2)\cr 1&0}$. 
It is equivalent to $(j(v_2))=(m_1)$, so $H_1(\widehat 
X_{\tau^D})=\Lambda/(j(v_2))=\Lambda/(m_1)$.
With the notation of Figure \ref{fig:2-tangles}, we have 
$\tau^N=(\tau_2)^D$. Hence, the formula for $\tau^N$ follows from the 
formula for $\tau^D$ and Proposition \ref{prop:move}.
\qed\end{pf}

\section{Generalizations}

\subsection{The category of m-colored tangles}

Fix throughout this section a positive integer $m$.
An {\em $m$-colored tangle} is an oriented tangle $\tau$ together with 
a map $c$ assigning
to each component $\tau_j$ of $\tau$ a {\em color} 
$c(j)\in\{1,\dots,m\}$. The composition of two
$m$-colored tangles is defined if and only if it is compatible with 
the coloring of each component.
Finally, we say that an $m$-colored tangle is an {\em oriented 
$m$-colored braid} if the underlying tangle is a braid.

More formally, $m$-colored tangles can be understood as morphisms of a 
category in the following way.
Consider two maps $\varphi\colon\{1,\dots,n\}\to\{\pm 1,\dots,\pm m\}$ 
and $\varphi'\colon\{1,\dots,n'\}\to\{\pm 1,\dots,\pm m\}$, where $n$ 
and $n'$ are non-negative integers. We will say that an $m$-colored 
tangle $(\tau,c)$ is a {\em $(\varphi,\varphi')$-tangle} if the 
following conditions hold:
\begin{itemize}
\item{$\tau$ is an $(\eps,\eps')$-tangle, where 
$\eps=\varphi/|\varphi|$ and $\eps'=\varphi'/|\varphi'|$;}
\item{if $x_i\in D^2\times\{0\}$ (resp. $x'_i\in D^2\times\{1\}$) is 
an endpoint of a component $\tau_j$ of $\tau$,
then $|\varphi(i)|=c(j)$ (resp. $|\varphi'(i)|=c(j)$).}
\end{itemize}
Two $(\varphi,\varphi')$-tangles are isotopic if they are isotopic as 
$(\eps,\eps')$-tangles under an isotopy that respects the color of
each component.
We denote by $T(\varphi,\varphi')$ the set of isotopy classes of 
$(\varphi,\varphi')$-tangles.
The composition of oriented tangles induces a composition
$T(\varphi,\varphi')\times T(\varphi',\varphi'')\to 
T(\varphi,\varphi'')$ for any $\varphi,\varphi'$ and $\varphi''$.

This allows us to define the {\em category of $m$-colored tangles} 
${\mathbf Tangles_m}$. Its objects are the maps
$\varphi\colon\{1,\dots,n\}\to\{\pm 1,\dots,\pm m\}$ with $n\ge 0$, 
and its morphisms are given by
$\Hom(\varphi,\varphi')=T(\varphi,\varphi')$. Clearly, oriented 
$m$-colored braids and oriented $m$-colored string links form 
categories ${\mathbf Braids_m}$ and ${\mathbf Strings_m}$ such that
\[
{\mathbf Braids_m}\subset{\mathbf Strings_m}\subset{\mathbf Tangles_m}.
\]

\subsection{The multivariable Lagrangian representation}

We now define a functor ${\mathscr 
F}_m\colon{\mathbf Tangles_m}\to{\mathbf Lagr_{\Lambda_m}}$, where 
$\Lambda_m$ denotes the ring $\Z[t_1^{\pm 1},\dots,t_m^{\pm 1}]$.
This construction generalizes the functor of Theorem \ref{thm:1}, 
which corresponds to the case $m=1$.
It also extends the works of Gassner for pure braids and Le Dimet for 
pure string links.

Consider an object of ${\mathbf Tangles_m}$, that is, a map
$\varphi\colon\{1,\dots,n\}\to\{\pm 1,\dots,\pm m\}$ with $n\ge 0$.
Set 
$\ell_\varphi=(\ell^{(1)}_\varphi,\dots,\ell^{(m)}_\varphi)\in\Z^m$,
where $\ell^{(j)}_\varphi=\sum_{\{i\mid\varphi(i)=\pm 
j\}}sign(\varphi(i))$ for $j=1,\dots,m$. Using the notation of Section 
3.2, we define
\[
D_\varphi=\cases{D^2\setminus\N(\{x_1,\dots,x_n\})& if 
$\ell_\varphi\neq(0,\dots,0)$;\cr
	S^2\setminus\N(\{x_1,\dots,x_n\})& if 
$\ell_\varphi=(0,\dots,0)$.}
\]
As in the case of oriented tangles, we endow $D_\varphi$ with the 
counterclockwise orientation, a base point $z$, and generators 
$e_1,\dots,e_n$ of $\pi_1(D_\varphi,z)$. Consider the homomorphism 
from $\pi_1(D_\varphi)$ to the free abelian group $G\cong\Z^m$ with 
basis $t_1,\dots,t_m$ given by
$e_i\mapsto t_{|\varphi(i)|}$. It defines a regular $G$-covering 
$\widehat{D}_\varphi\to D_\varphi$, so the homology 
$H_1(\widehat{D}_\varphi)$ is a module over $\Z G=\Lambda_m$.
Finally, let $\omega_\varphi\colon H_1(\widehat{D}_\varphi)\times 
H_1(\widehat{D}_\varphi)\to\Lambda_m$ be the skew-hermitian pairing 
given by
\[
\omega_\varphi(x,y)=\sum_{g\in G}<gx,y> g^{-1},
\]
where $<\;,\;>\colon H_1(\widehat{D}_\varphi)\times 
H_1(\widehat{D}_\varphi)\to\Z$ is the intersection form induced by the 
orientation of
$D_\varphi$ lifted to $\widehat{D}_\varphi$.

Consider now a $(\varphi,\varphi')$-tangle $(\tau,c)$. Note that 
$\ell_\varphi=\ell_{\varphi'}$. Let $X_\tau$ be the compact manifold
\[
X_\tau=\cases{(D^2\times [0,1])\setminus\N(\tau)& if 
$\ell_\varphi\neq(0,\dots,0)$;\cr
	(S^2\times [0,1])\setminus\N(\tau)& if 
$\ell_\varphi=(0,\dots,0)$,} 
\]
oriented so that the induced orientation on $\partial X_\tau$ extends 
the orientation on $(-D_\varphi)\sqcup D_{\varphi'}$.
We know from Section 3.3 that $H_1(X_\tau)=\bigoplus_{j=1}^\mu \Z m_j$ 
if $\ell_\varphi\neq(0,\dots,0)$, and
$H_1(X_\tau)=\bigoplus_{j=1}^\mu \Z m_j/\sum_{i=1}^n 
sign(\varphi(i))e_i$ otherwise.
Hence, the coloring of $\tau$ defines a homomorphism $H_1(X_\tau)\to 
G$,
$m_j\mapsto t_{c(j)}$ which induces a homomorphism $\pi_1(X_\tau)\to 
G$ extending the homomorphisms $\pi_1(D_\varphi)\to G$ and
$\pi_1(D_{\varphi'})\to G$. It gives a $G$-covering 
$\widehat{X}_\tau\to X_\tau$.

Consider the inclusion homomorphisms $i_\tau\colon 
H_1(\widehat{D}_\varphi)\to H_1(\widehat{X}_\tau)$ and
$i'_\tau\colon H_1(\widehat{D}_{\varphi'})\to H_1(\widehat{X}_\tau)$. 
Denote by $j_\tau$ the homomorphism
$H_1(\widehat{D}_\varphi)\oplus H_1(\widehat{D}_{\varphi'})\to 
H_1(\widehat{X}_\tau)$ given by $j_\tau(x,x')=i'_\tau(x')-i_\tau(x)$. 
Set
\[
{\mathscr F}_m(\tau)=\ker(j_\tau)\subset 
H_1(\widehat{D}_\varphi)\oplus H_1(\widehat{D}_{\varphi'}). 
\]
\begin{thm}\label{thm:multi}
Let ${\mathscr F}_m$ assign to each map 
$\varphi\colon\{1,\dots,n\}\to\{\pm 1,\dots,\pm m\}$ the pair
$(H_1(\widehat D_\varphi),\omega_\varphi)$ and to each $\tau\in 
T(\varphi,\varphi')$ the submodule ${\mathscr F}_m(\tau)$ of
$H_1(\widehat D_\varphi)\oplus H_1(\widehat D_{\varphi'})$. Then, 
${\mathscr F}_m$ is a functor 
${\mathbf Tangles_m}\to{\mathbf Lagr_{\Lambda_m}}$ which fits in the 
diagrams
\[
\begin{CD}
{\mathbf Tangles_m}@>{{\mathscr F}_m}>>{\mathbf Lagr_{\Lambda_m}}\\
@AAA @AAA\\
{\mathbf Braids_m}@>>>{\mathbf U_{\Lambda_m}}
\end{CD}
\qquad and\qquad
\begin{CD}
{\mathbf Tangles_m}@>{{\mathscr J}\circ{\mathscr 
F}_m}>>{\mathbf \overline{Lagr}_{\Lambda_m}}\\
@AAA @AAA\\
{\mathbf Strings_m}@>>>{\mathbf U^0_{\Lambda_m}}
\end{CD},
\]
where the vertical arrows denote embeddings of categories.
\end{thm}
\begin{pf}
Lemmas \ref{lemma:non_deg}, \ref{lemma:Lag'}, \ref{lemma:comp}, 
Proposition \ref{prop:braids}, Proposition \ref{prop:stringlinks} and 
their proofs extend to our setting with obvious changes.
The only `topological' facts required are the following:
\begin{enumerate}
\item{$H_1(\partial\widehat{D}_\varphi)=0$,}
\item{the $\Lambda_m$-module $H_1(\widehat{D}_\varphi)$ is 
torsion-free,
}
\item{$H_1(\partial\widehat{X}_\tau,\widehat{D}_\varphi\sqcup
\widehat{D}_{\varphi'})$ is a torsion $\Lambda_m$-module.}
\end{enumerate}
The definition of $\widehat{D}_\varphi$ easily implies that 
$\partial\widehat{D}_\varphi$ consists of copies of $\R$, so the first 
claim is checked. Since
$D_\varphi$ has the homotopy type of a $1$-dimensional $CW$-complex 
$Y_\varphi$, the $\Lambda_m$-module $H_1(\widehat{D}_\varphi)=
H_1(\widehat{Y}_\varphi)=Z_1(\widehat{Y}_\varphi)$ is a submodule of 
the free $\Lambda_m$-module $C_1(\widehat{Y}_\varphi)$. Therefore,
$H_1(\widehat{D}_\varphi)$ is torsion-free. Finally, the third claim 
follows easily from the definitions and the
excision theorem.
\qed\end{pf}

\subsection{High-dimensional Lagrangian representations}

The Lagrangian representation of Theorem \ref{thm:1} can be 
generalized in another direction by considering high-dimensional 
manifolds. We conclude the paper with a brief sketch of this 
construction.

Fix throughout this section an integer $n\ge 1$. In the sequel, all 
the manifolds are assumed piecewise linear, compact and oriented. 
Consider a homology $2n$-sphere $D$. To this manifold, we associate a 
category
${\mathscr C}_D$ as follows. Its objects are codimension-$2$ 
submanifolds $M$ of $D$ such that $H_n(M)=0$. The morphisms between 
$M\subset D$ and $M'\subset D$ are given by properly embedded 
codimension-$2$ submanifolds $T$ of $D\times [0,1]$ such that
the oriented boundary $\partial T$ of $T$ satisfies $\partial T\cap 
(D\times\{0\})=-M$ and $\partial T\cap (D\times\{1\})=M'$, where $-M$ 
denotes $M$ with the opposite orientation. The composition is defined 
in the obvious way.

If $D_M$ is the complement of an open tubular neighborhood of $M$ in 
$D$, we easily check that $H_1(D_M)\cong H_0(M)$. Therefore, the 
epimorphism $H_0(M)\to\Z$ which sends every generator to $1$ 
determines a $\Z$-covering $\widehat D_M\to D_M$. The lift of the 
orientation of $D_M$ to $\widehat D_M$ defines a $\Z$-bilinear 
intersection form on $H_n(\widehat D_M)$. This gives a 
$\Lambda$-sesquilinear form on $H_n(\widehat D_M)$, which in turn 
induces a $\Lambda$-sesquilinear form $\omega_M$ on $BH_n(\widehat 
D_M)$, where $BH=H/Tors_\Lambda H$ for a $\Lambda$-module $H$. (Note 
that $\omega_M$ is skew-hermitian if $n$ is odd, and Hermitian if $n$ 
is even.) Using the fact that $H_n(M)=0$, the proof of Lemma 
\ref{lemma:non_deg} can be applied to this setting, showing that 
$\omega_M$ is non-degenerate. Let ${\mathscr F}_D(M)$ denote the 
$\Lambda$-module $BH_n(\widehat D_M)$ endowed with the non-degenerate 
$\Lambda$-sesquilinear form $\omega_M$.

Given a codimension-$2$ submanifold $T$ of $D\times[0,1]$, denote by 
$X_T$ the complement of an open tubular neighborhood of $T$ in 
$D\times[0,1]$. Since $H_1(X_T)\cong H_0(T)$, we have a $\Z$-covering 
$\widehat X_T\to X_T$ given by the homomorphism $H_0(T)\to\Z$ which 
sends every generator to $1$. There are obvious inclusions $\widehat 
D_M\subset\widehat X_T$ and $\widehat D_{M'}\subset\widehat X_T$ which 
induce homomorphisms $i$ and $i'$ in $n$-dimensional homology. Let 
$j\colon H_n(\widehat D_{M})\oplus H_n(\widehat D_{M'})\to 
H_n(\widehat X_{T})$ be the homomorphism given by 
$j(x,x')=i'(x')-i(x)$. It induces a homomorphism
\[
BH_n(\widehat D_M)\oplus BH_n(\widehat 
D_{M'})\stackrel{j_T}{\longrightarrow}BH_n(\widehat X_{T}).
\]
Set ${\mathscr F}_D(T)=\ker(j_T)$. The proof of Lemma \ref{lemma:Lag'} 
can be applied to check that ${\mathscr F}_D(T)$ is a Lagrangian 
submodule of $(-BH_n(\widehat D_M))\oplus BH_n(\widehat D_{M'})$. 
Lemma \ref{lemma:comp} can also be adapted to our setting to show that 
${\mathscr F}_D(T_2\circ T_1)={\mathscr F}_D(T_2)\circ{\mathscr 
F}_D(T_1)$. Therefore, ${\mathscr F}_D$ is a functor from
${\mathscr C}_D$ to the Lagrangian category 
${\mathbf \overline{Lagr}_\Lambda}$ amended as follows: the 
non-degenerate form is Hermitian if $n$ is even, skew-hermitian if $n$ 
is odd.

\begin{ack}
A part of this paper was done while the authors visited the Department 
of Mathematics 
of the Aarhus University whose hospitality the authors thankfully 
acknowledge. The first author also wishes to thank the
Institut de Recherche Math\'ematique Avanc\'ee (Strasbourg) and the
Institut de Math\'ematiques de Bourgogne (Dijon) for hospitality.
\end{ack}

\end{document}